\documentclass[11pt]{amsart}

\usepackage[T1]{fontenc}
\usepackage[utf8]{inputenc}
\usepackage[english]{babel}
\usepackage{amscd,amsmath,amsthm,amssymb}
\usepackage{mathtools}
\usepackage{mathrsfs}
\usepackage{comment}
\usepackage{xcolor}
\usepackage{array}
\usepackage{caption}
\usepackage{mathabx}
\usepackage{standalone}
\usepackage{scalerel}
\usepackage{csquotes}

\usepackage{enumitem}
\setlist[enumerate]{label=\textnormal{(\arabic*)}}
\newcommand{\defaultRoman}{{\textnormal{(\textit{\roman*})}}}

\definecolor{cadmiumgreen}{rgb}{0.0, 0.42, 0.24}
\usepackage[
  colorlinks,
  citecolor=cadmiumgreen,
  pdfstartview ={FitV},
]{hyperref}
\hypersetup{
  pdfauthor={Omid Amini, Matthieu Piquerez},
  pdftitle={Tropical Clemens-Schmid and existence of tropical cycles}}

\usepackage[
  alphabetic,
  msc-links,
  nobysame,
  lite,
]{amsrefs}

\usepackage{tikz, tikz-cd}
\usetikzlibrary{calc, backgrounds, quotes, arrows.meta, positioning, decorations.pathmorphing}

\usepackage[left=3cm,top=3.5cm,right=3cm]{geometry}
\allowdisplaybreaks[1]

\newtheorem{thm}{Theorem}[section]
\newtheorem{lemma}[thm]{Lemma}
\newtheorem{prop}[thm]{Proposition}

\newtheorem{conj}[thm]{Conjecture}

\theoremstyle{definition}

\newtheorem{remm}[thm]{Remark}
\newenvironment{remark}
  {\pushQED{\qed}\remm}
  {\popQED\endremm}

\numberwithin{equation}{section}

\newcommand{\cf}[1]{cf.}
\newcommand{\ie}{i.e.}

\renewcommand{\~}{\widetilde}
\newcommand{\Q}{\mathbb{Q}}
\newcommand{\Z}{\mathbb{Z}}
\newcommand{\R}{\mathbb{R}}

\newcommand{\simto}{\xrightarrow{\raisebox{-3pt}[0pt][0pt]{\small$\hspace{-1pt}\sim$}}}
\newcommand{\longsimto}{\xrightarrow{\ \raisebox{-3pt}[0pt][0pt]{\small$\hspace{-1pt}\sim$\ }}}
\newcommand{\bul}{\bullet} 

\renewcommand{\emptyset}{\varnothing}
\newcommand{\rquot}[2]{#1\big/#2}
\newcommand{\rest}[1]{\raisebox{-1pt}{$\vert$}_{#1}}
\renewcommand{\pmod}[1]{\ (\mathrm{mod}\ #1)}

\newcommand{\ccdot}{\,\cdot\,}
\newcommand{\rdot}{\cdot\,}

\let\oldchi\chi
  \newcommand{\raisechi}[2]{\raisebox{.4ex}{$#1#2$}}
  \renewcommand{\chi}{{\mathpalette\raisechi\oldchi}}

\makeatletter
\let\oldsum\sum
\renewcommand{\sum}{\@ifnextchar_\@mysum\oldsum}
\def\@mysum_#1{\oldsum_{\substack{#1}}}
\let\oldbigoplus\bigoplus
\renewcommand{\bigoplus}{\@ifnextchar_\@mybigoplus\oldbigoplus}
\def\@mybigoplus_#1{\oldbigoplus_{\substack{#1}}}
\let\oldprod\prod
\renewcommand{\prod}{\@ifnextchar_\@myprod\oldprod}
\def\@myprod_#1{\oldprod_{\substack{#1}}}
\makeatother

\let\oldbigwedge\bigwedge
\renewcommand{\bigwedge}{{\textstyle\oldbigwedge\!}}

\newcommand{\cdvdots}{\raisebox{0pt}[14pt][4pt]{\makebox[10pt]{$\vdots$}}}

\let\Im\relax
  \DeclareMathOperator{\Im}{Im} 
\DeclareMathOperator{\coker}{coker}
\DeclareMathOperator{\Hom}{Hom} 
\DeclareMathOperator{\sed}{sed} 
\DeclareMathOperator{\gys}{Gys} 
\DeclareMathOperator{\id}{id} 
\DeclareMathOperator{\class}{cl} 
\DeclareMathOperator{\Cone}{Cone} 
\DeclareMathOperator{\conv}{conv}
\DeclareMathOperator{\sign}{sign}

\newcommand{\trop}{{\scaleto{\mathrm{trop}}{5pt}}} 
\newcommand{\Dolb}{{\scaleto{\mathrm{Dolb}}{5pt}}} 
\newcommand{\RpMod}{\mathcal M} 
\newcommand{\eR}{\mathbf R} 
\newcommand{\Cint}{\mathring{C}} 
\newcommand{\e}{{\mathfrak e}} 
\newcommand{\f}{{\textnormal{\textsl{\textsf f}}}} 
\newcommand{\II}{\mathcal I} 
\renewcommand{\P}{\mathbb P} 
\newcommand{\TP}{\mathbb{TP}} 
\newcommand{\TT}{\mathrm{T}} 
\newcommand{\dual}{\star} 
\let\i\relax
  \newcommand{\i}{{\mathop{}\mathrm{i}}} 
\renewcommand{\d}{{\rm d}} 
\newcommand{\SF}{\textrm{\bf F}} 
\newcommand{\Cycle}{\mathscr C} 
\newcommand{\MW}{\mathrm{MW}} 
\newcommand{\nvect}{\mathfrak n} 
\newcommand{\x}{\textsc{x}} 

\newcommand{\HL}{\textrm{HL}} 
\newcommand{\PD}{\textrm{PD}} 

\newcommand{\st}{\bigm|} 
\newcommand{\Bigst}{\Bigm|} 
\newcommand{\abs}[1]{\lvert #1\rvert} 
\newcommand{\comp}[1]{\overline{#1}} 
\newcommand{\sminfty}{{\scaleobj{.7}{\infty}}} 

\newcommand{\X}{\mathfrak X}
\newcommand{\Y}{\mathscr Y}

\newcommand{\dimsaux}[2]{\raisebox{.2ex}{\scalebox{1}[.8]{$#1\lvert$}}#2\raisebox{.2ex}{\scalebox{1}[.8]{$#1\rvert$}}}
  \newcommand{\dims}[1]{\mathpalette\dimsaux{#1}}
\newcommand{\suppaux}[2]{\scalebox{1}[1.4]{$#1\lvert$}#2\scalebox{1}[1.4]{$#1\rvert$}}
  \newcommand{\supp}[1]{\mathpalette\suppaux{#1}}
\newcommand{\conezero}{{\underline0}} 

\newcommand{\subface}{\prec}
\newcommand{\ssubface}{\mathbin{\mathchoice
  {\subface\!\!\!\cdot}%
  {\subface\!\!\!\cdot}%
  {\subface\!\cdot}%
  {\subface\!\cdot}%
}} 
\newcommand{\supface}{\succ}
\newcommand{\ssupface}{\mathbin{\mathchoice
  {\cdot\!\!\!\supface}%
  {\cdot\!\!\!\supface}%
  {\cdot\!\supface}%
  {\cdot\!\supface}%
}}

\newcommand{\Ma}{\mathfrak M} 
\newcommand{\Fl}{\mathscr{F}}

\newcommand{\ST}{\textnormal{\textsf{S\!T}}} 
\renewcommand{\AA}{\textnormal{\textsf{A}}}

\newcommand{\fX}{\mathscr X} 
\newcommand{\ufX}{\~{\mathscr X^*}} 
\newcommand{\disk}{\triangle} 
\newcommand{\pdisk}{\disk^{\!*}} 
\newcommand{\udisk}{\~\pdisk} 

\newcommand{\s}{{\rm s}} 
\newcommand{\rel}{{\rm rel}} 

\begin{document}
\title[Tropical Clemens-Schmid and existence of tropical cycles]{Tropical Clemens-Schmid sequence and existence of tropical cycles with a given cohomology class}

\author{Omid Amini}
\address{CNRS - CMLS, \'Ecole Polytechnique}
\email{\href{omid.amini@polytechnique.edu}{omid.amini@polytechnique.edu}}

\author{Matthieu Piquerez}
\address{CMLS, \'Ecole Polytechnique}
\email{\href{matthieu.piquerez@polytechnique.edu}{matthieu.piquerez@polytechnique.edu}}

\date{December 24, 2020}

\begin{abstract} This is a sequel to our work in tropical Hodge theory. Our aim here is to prove a tropical analogue of the Clemens-Schmid exact sequence in asymptotic Hodge theory. As an application of this result, we prove the tropical Hodge conjecture for smooth projective tropical varieties which are rationally triangulable.
 This provides a partial answer to a question of Kontsevich who suggested the validity of the tropical Hodge conjecture could be used as a test for the validity of the Hodge conjecture.
\end{abstract}
\maketitle

\setcounter{tocdepth}{1}

\tableofcontents

\section{Introduction and statement of the main theorem}

\subsection{Tropical Hodge conjecture} In this paper we prove the following theorem.

\begin{thm}[Hodge conjecture for rationally triangulable smooth projective tropical varieties]\label{thm:HC} Let $\X$ be a smooth projective tropical variety. Assume that $\X$ is rationally triangulable. The locus of Hodge classes in $H^{p,p}_\trop(\X, \Q)$ associated to codimension $p$ tropical cycles in $\X$ coincides with the kernel of the tropical monodromy map $N\colon H^{p,p}_\trop(\X, \Q) \to H^{p-1, p+1}_\trop(\X,\Q)$.
\end{thm}
The notions of tropical smoothness and rational triangulability, and the definition of the tropical cohomology groups will be recalled in Section~\ref{sec:preliminaries}.

\medskip

In~\cite{Zha20}, Zharkov explains a suggestion of Kontsevitch on how to test the validity of the Hodge conjecture via integral affine manifolds by specialization: if the tropical Hodge conjecture turned out to be false for some tropical limit of abelian varieties, for example, then this would imply that the classical conjecture would be false as well. On the other hand, the validity of the tropical Hodge conjecture for general tropical varieties might be regarded as an evidence for the validity of the Hodge conjecture. Our Theorem~\ref{thm:HC} above goes in the direction of this suggestion. While it does not answer the question of Kontsevich in its original form, which concerns tropical abelian varieties which are in general not rationally triangulable, it goes beyond the case of affine manifolds.

A stronger form of the above theorem in codimension one, the tropical analogue of Lefschetz $(1,1)$-theorem, with integral coefficients and without the triangulability assumption, was proved by Jell-Rau-Shaw~\cite{JRS19}. The general form of the tropical Hodge conjecture, without the triangulability assumption, is as follows.

\begin{conj}[Tropical Hodge conjecture]
Let $\X$ be a smooth projective tropical variety. The locus of Hodge classes in $H^{p,p}_\trop(\X,\Q)$ generated by classes of codimension $p$ tropical cycles in $\X$ coincides with the kernel of the tropical monodromy operator $N\colon H^{p,p}_\trop(\X, \Q) \to H^{p-1, p+1}_\trop(\X,\R)$.
\end{conj}
In the case where $\X$ is rationally triangulable, the monodromy operator is rational, and the statement above is the content of Theorem~\ref{thm:HC}. In general, $N$ is only defined with real coefficients.

\smallskip

It is well-known that the Hodge conjecture implies the Grothendieck's standard conjecture of type D, that the numerical and homological equivalence on algebraic cycles coincide, for varieties over a field of characteristic zero. We prove the tropical analogue of this standard conjecture.
\begin{thm}\label{thm:standard} Let $\X$ be a smooth projective tropical variety which is rationally triangulable. The numerical and homological equivalences on tropical cycles coincide.
\end{thm}
This answers a question of Gross and Shokrieh~\cite{GS19-taut} for rationally triangulable smooth projective tropical varieties.

\subsection{Tropical Clemens-Schmid} \label{sec:classical}

In order to prove the above results, we prove a tropical analogue of the Clemens-Schmid exact sequence in asymptotic Hodge theory~\cites{Cle77,Sch73} and its extension to algebraic cycles by Bloch-Gilet-Soul\'e~\cite{BGS95-degenerate-cycles}. This will be based on results we proved in our paper~\cite{AP-tht}, which we will recall in Sections~\ref{sec:preliminaries} and~\ref{sec:steenbrink}.

We start by recalling the classical Clemens-Schmid sequence. Let $\fX^*$ be a projective family of smooth complex varieties over the punctured disk $\pdisk$. Passing to a finite \'etale cover of the punctured disk if necessary, we can complete $\fX^*$ to a regular semistable family $\fX$ over the disk $\disk$. This gives a special fiber $\fX_0$ over $0$, whose addition as relative boundary results in a relative compactification $\fX$ of $\fX^*$ over the disk $\disk$.

Let $\udisk \to \pdisk$ be the universal cover of the punctured disk $\pdisk$ and denote by $\ufX$ the family of complex varieties over $\udisk$ obtained by pulling back $\fX^*$ over $\udisk$.

There is a monodromy operator $T\colon \ufX\to\ufX$ associated to the generator of the fundamental group of the punctured disk $\pdisk$, which induces an automorphism on each fiber of the original family. The induced operator $T$ on the cohomology of $\ufX$ is unipotent and leads to the (logarithmic) monodromy operator $N := -\log(T)\colon H^\bul(\ufX) \to H^\bul(\ufX)$. The operator $N$ is nilpotent and the corresponding Jacobson-Morosov filtration together with an appropriate Hodge filtration endows the cohomology $H^\bul(\ufX)$ of $\ufX$ with a mixed Hodge structure called the \emph{limit mixed Hodge structure} of the family~\cite{Sch73}. This limit mixed Hodge structure can be computed algebraically thanks to the Steenbrink spectral sequence~\cite{Ste76}.

\medskip

Since the family $\fX$ retracts by deformation to $\fX_0$, we get an isomorphism of the cohomology $H^\bul(\fX)$ with that of $\fX_0$. For this reason, these cohomology groups might be named the \emph{surviving cohomology} of the family.
Each surviving cohomology group comes with a mixed Hodge structure which can be defined using the Deligne spectral sequence~\cite{Del-hodge2}.

\medskip

The open inclusion $\fX^* \hookrightarrow \fX$ leads to the definition of the cohomology groups $H^\bul(\fX, \fX^*)$, which can be named the \emph{relative cohomology} of the family. Again, Deligne's theory endows these cohomology groups with mixed Hodge structures.

\medskip

The Clemens-Schmid exact sequence~\cite{Cle77} is an exact sequence which establishes a link between the above mixed Hodge structures, on various degrees. This is the following long exact sequence of mixed Hodge structures
\[ \cdots \to H^k(\fX) \to H^k(\ufX) \xrightarrow{N} H^k(\ufX) \to H^{k+2}(\fX,\fX^*) \to H^{k+2}(\fX) \to \cdots, \]
where the morphisms in the sequence are of specific degrees that we do not precise here.

\medskip

We also have two distinguished triangles in the derived category:
\[ \begin{tikzcd}[column sep=tiny]
H^\bul(\ufX) \ar[rr, "N"] && H^\bul(\ufX) \ar[dl, "+1"] \\
 & \ar[ul] H^\bul(\fX^*) &
\end{tikzcd} \qquad \begin{tikzcd}[column sep=tiny]
H^\bul(\fX) \ar[rr] && H^\bul(\fX^*) \ar[dl, "+1"] \\
 & \ar[ul] H^\bul(\fX, \fX^*) &
\end{tikzcd}
\]

\medskip

In this paper, we define the tropical analogue of the \emph{surviving} and \emph{relative} cohomology groups for a smooth projective tropical variety $\X$.

The definition of these cohomology groups is dependent on the choice of a unimodular triangulation $X$ on $\X$ (which exists after changing the underlying lattice by a rational multiple, by the rational triangulability assumption). For each pair of non-negative integers $p,q$, we will define the bigraded surviving and relative cohomology groups $H^{p,q}_{s}(X, \Q)$ and $H^{p,q}_\rel(X, \Q)$, respectively, and set for any non-negative integer $k$,
\[H^k_s(X) := \bigoplus_{p+q =k} H^{p,q}_{s}(X, \Q) \qquad \textrm{and} \qquad H^k_\rel(X) := \bigoplus_{p+q =k} H^{p,q}_{\rel}(X, \Q).\]
These cohomology groups come with canonical maps $H^k_s(X) \to H^k(\X)$ and $H^k_\rel(X) \to H^k_s(X)$. Here, we set
\[H^k(\X) := \bigoplus_{p+q=k} H^{p,q}_\trop(\X, \Q)\]
and note that it does not depend on the choice of the triangulation.

The above cohomology groups fit together and give the following long exact sequence.

\begin{thm}[Tropical Clemens-Schmid exact sequence] \label{thm:CS} We have the following exact sequence
\[ \dots \to H^k_s(X) \to H^k(\X) \xrightarrow{N} H^k(\X) \to H^{k+2}_\rel(X) \to H^{k+2}_s(X) \to H^{k+2}(\X) \xrightarrow{N} H^{k+2}(\X) \to \cdots
\]
where the morphism $N\colon H^k(\X) \to H^{k}(\X)$ is given by the sum of the tropical monodromy maps $N \colon H^{p,q}_{\trop}(\X, \Q) \to H^{p-1, q+1}_\trop(\X, \Q)$.
\end{thm}

\subsection{Explicit description of the tropical cycle associated to a Hodge class} Let $\X$ be a smooth projective tropical variety. Suppose that $\X$ is rationally triangulable. Let $\alpha$ be a tropical Hodge class, that is $\alpha \in H^{p,p}(\X, \Q) \cap \ker(N)$, for some non-negative integer $p$. We fix a unimodular triangulation $X$ of $\X$ (which again exists after changing the underlying lattice by a rational multiple). Using the Clemens-Schmid exact sequence, we construct an explicit tropical cycle $C$ with cohomology class equal to $\alpha$. The tropical cycle $C$ is defined by a \emph{Minkowski weight} on the $p$-dimensional cells of the triangulation $X$, which is itself obtained by gluing local Minkowski weights around vertices of the triangulation. The existence of local Minkowski weights is a consequence of the tropical Clemens-Schmid exact sequence. The fact that the image of $C$ by the tropical cycle class map coincides with $\alpha$ is a consequence of our proof of the Steenbrink-Tropical comparison theorem~\cite{AP-tht}.

\subsection{Organization of the paper} In Section~\ref{sec:preliminaries}, we recall basis results and definitions concerning tropical varieties. In Section~\ref{sec:minkowski_weight}, we define Minkowski weights in the local and global setting. In Section~\ref{sec:local_HC}, we discuss the tropical Hodge conjecture in the local setting. Section~\ref{sec:steenbrink} recalls some results concerning the tropical Steenbrink sequence and the comparison theorem from~\cite{AP-tht} by making them more explicit for some specific primitive parts which will be of later use in the study of the tropical cycle class map. In Section~\ref{sec:surviving_relative} we introduce the tropical surviving and relative cohomology groups, relative to a triangulation of the tropical variety. In Section~\ref{sec:CS_sequence}, we define the tropical Clemens-Schmid sequence and prove Theorem~\ref{thm:CS}. Section~\ref{sec:proof_HC} contains the proof of Theorem~\ref{thm:HC}. Finally, Theorem~\ref{thm:standard} is proved in Section~\ref{sec:standard}.

\section{Preliminaries on tropical varieties} \label{sec:preliminaries}

The aim of this section is to provide a brief account of tropical geometry and the necessary background on the terminology used in the statement of our main theorems. We refer to~\cite{AP-tht} for more details.

\subsection{Polyhedral complexes and star fans} \label{sec:recol} We start by fixing the polyhedral geometry terminology and introducing (extended) polyhedral complexes.

\medskip

A polyhedron $P$ in a real vector space $V\simeq \R^n$ is by definition a non-empty intersection of a finite number of affine half-spaces. We consider only polyhedra which are strongly convex in the sense that they do not contain any affine line. The tangent space of $P$ that we denote by $\TT(P)$ or by $\TT P$ is the linear subspace of $V$ spanned by the differences $x-y$ for pairs of elements $x,y$ of $P$. We denote by $\dims{P}$ the dimension of $P$. A face of a polyhedron $P$ is either $P$ itself or a nom-empty intersection of $P$ with an affine hyperplane $H$ provided that $P$ lies in one of the two half-spaces delimited by $H$. For two polyhedra $\gamma$ and $\delta$, we write $\gamma\subface\delta$ to indicate that $\gamma$ is a face of $\delta$. We use the notation $\gamma\ssubface\delta$ if moreover $\gamma$ is of codimension one in $\delta$. The partial order $\subface$ puts a lattice structure on the set of faces of a polyhedron. A face of dimension zero in $P$ is called a vertex and a face of dimension one is called an edge.

By a cone we mean a strongly convex polyhedron with a unique vertex which is the origin of $\R^n$. A polytope is a compact polyhedron.

Assume now that the real vector space $V \simeq \R^n$ comes with a lattice $N \simeq \Z^n$ of full rank, so that $N_\R := N\otimes_\Z \R = V$. In such a case, we say a polyhedron $P$ in $V$ is rational if all the half-spaces which are used to define $P$ can be defined in $N_\Q := N \otimes_\Z \Q \simeq \Q^n$. For a rational polyhedron $P$, we denote by $N_P:=N\cap \TT(P)$ the full-rank lattice of the tangent space $\TT(P)$.

We say $P$ is integral with respect to $N$ if it is rational and its vertices are all in $N$. If the lattice $N$ is understood from the context, we omit to mention it.

A polyhedral complex $Y$ in a real vector space $V$ is a finite non-empty collection of polyhedra in $V$ called faces of $Y$ such that for any pair of faces $\delta$, $\delta'$ in $Y$, the following two properties hold.
\begin{enumerate}[label=\defaultRoman]
\item \label{enum:subface} Any face of $\delta$ is contained in $Y$.
\item The intersection $\delta\cap\delta'$ is either empty or is a face of both $\delta$ and $\delta'$.
\end{enumerate}

A fan is a polyhedral complex $\Sigma$ which has a unique vertex the origin of $V$. In this case, all the faces of $\Sigma$ are cones. A cone of dimension one is called a ray.

For a polyhedral complex $Y$ and for a non-negative integer $k$, we denote by $Y_{(k)}$ the $k$-skeleton of $Y$ which consists of all the faces of $Y$ of dimension at most $k$ and by $Y_k$ the set of $k$-dimensional faces of $Y$. The dimension $d$ of $Y$ is the maximum of the dimension of its faces. The elements of $Y_d$ are called facets. We say that $Y$ is of pure dimension $d$ if every face of $Y$ is included in some facet of $Y$.

The support of a polyhedral complex $Y$ denoted by $\supp{Y}$ is the union of all the faces of $Y$ in $V$.

A polyhedral complex $Y$ is called rational, integral, unimodular with respect to the lattice $N$ if all faces of $Y$ are rational, integral, or unimodular with respect to $N$, respectively. We already defined the terminology rational and integral, we recall what unimodular means here.

First, recall that for two polyhedra $A$ and $B$ in $V$, the Minkowski sum $A+B$ is by definition the polyhedron
\[A+B := \bigl\{x+ y \st x\in A, y \in B\bigr\}.\]
Every polyhedron $P$ can be written as the Minkowski sum $Q+\sigma$ of a polytope $Q$ and a cone $\sigma$ by Minkowski-Weyl theorem. The cone $\sigma$ in the above decomposition is in addition unique and we will denote it by $P_\infty$; for $Q$ we can choose the convex-hull of the vertices of $P$. In particular, this implies that we can write \[ P=\conv(v_0,\dots,v_k)+\sum_{i=1}^l\R_+u_i\]
for points $v_0,\dots,v_k$ and vectors $u_1,\dots,u_l$ in $V$, where $\R_+$ denotes the space of non-negative real numbers. A polyhedron $P$ is called unimodular if it is integral with respect to $N$ and moreover, the points $v_i$ and the vectors $u_j$ can be chosen in $N$ in such a way that the collection of vectors $(v_1-v_0,\dots,v_k-v_0,u_1,\dots,u_l)$ form a basis of the lattice $N_P = N\cap\TT(P)$. In this case, the points $v_0,\dots, v_k$ coincide with the vertices of $P$ so that defining $P_\f:=\conv(v_0,\dots, v_k)$, we get the decomposition $P = P_\f + P_\infty$. In this decomposition, any point $x$ can be written in a unique way as the sum of two points $x_\f\in P_\f$ and $x_\infty \in P_\infty$. If we relax the above condition and only ask the collection of vectors $(v_1-v_0,\dots,v_k-v_0,u_1,\dots,u_l)$ to be independent we say the polyhedron $P$ is simplicial.

\medskip

For a free $\Z$-module $N$ of finite rank, we denote by $M=N^{\vee}$ its dual. The real vectors spaces corresponding to $N$ and $M$ are denoted by $N_\R$ and $M_\R$, respectively, and we have $M_\R =N_\R^\dual$. Note that we use the following convention in this paper: every time we work with a vector space, we use ${}^\dual$ for the dual vector space, and if we deal with a free $\Z$-module or a cone, we use instead ${}^\vee$ for the dual $\Z$-module and cone, respectively. Recall that for a rational polyhedron $\delta$ in $N_\R$, we use the notation $N_\delta$ to denote the lattice $N \cap \TT\delta$. We thus get the equality $N_{\delta, \R} = \TT\delta$.
Furthermore, we define the normal vector space of $\delta$ which we denote $N^{\delta}_{\R}$ by taking the quotient $N^{\delta}_{\R} := \rquot{N_\R}{N_{\delta, \R}}$. It comes with a full rank lattice $N^\delta = \rquot{N}{N_\delta}$.

\medskip

Let $X$ be a rational polyhedral complex in $N_\R$. For any face $\delta$ of $X$, one can choose a basis of $N_\delta$. The exterior product of the elements of this basis gives a generator of $\bigwedge^{\dims\delta}N_\delta$. This element is well-defined (that is independent of the choice of the basis) up to a sign. The choice of such an element for each face defines an orientation on $X$. In this article, we assume that every polyhedral complex is endowed with a fixed orientation. We denote by $\nvect_\delta$ the generator of $\bigwedge^{\dims\delta}N_\delta$ corresponding to this orientation and we call it the canonical unit multivector associated to $\delta$. Moreover, we denote by $\nvect_\delta^\dual$ its dual that we call the canonical $\dims\delta$-form associated to $\delta$. Note that $\nvect_\delta^\dual$ lives in the dual space $\bigwedge^{\dims\delta}N_\delta^\vee$.

\smallskip
Let $\gamma\ssubface\delta$ be  a pair of faces with $\gamma$ of codimension one in $\delta$. There is a unique generator $u_{\delta/\gamma}$ of $\rquot{N_\delta}{N_\gamma} \simeq \Z$ that lives in the part corresponding to $\delta$. We denote by $\nu^\dual_{\delta/\gamma}$ the linear form on $\TT\delta$ induced by the one on $\rquot{\TT\delta}{\TT\gamma}\simeq\R$ which takes value one at $u_{\delta/\gamma}$. The chosen orientation on $X$ induces a sign function on pair of faces $\gamma\ssubface\delta$ defined by the equation
\[ \nvect^\dual_\delta = \sign(\gamma, \delta) \nvect^\dual_\gamma \wedge \nu^\dual_{\delta/\gamma}, \]
where $\sign(\gamma, \delta) \in \{-1,+1\}$.

\medskip

We will follow the terminology introduced in our paper~\cite{AP-tht}: we will use a face $\delta$ as a subscript for subspaces of the ambient space or lattice or elements associated to these subspaces. And we use a face $\delta$ as a superscript to denote the quotient by $N_\delta$ or $N_{\delta, \R}$ of the ambient lattice or space, or to denote elements related to this quotient.

For a fan $\Sigma$ and a cone $\tau \in \Sigma$, the star fan of $\tau$ denoted by $\Sigma^\tau$ is defined by
\[\Sigma^\tau := \Bigl\{\, \pi_\tau (\sigma) \Bigst \sigma \supface \tau \textrm{ is a cone in $\Sigma$} \, \Bigr\}. \]
Here $\pi_\tau\colon N_{\R} \to N^\tau_\R$ is the projection map to the quotient. Our star fans are sometimes called transversal fans in other places, e.g., in~\cites{Kar04, BBFK02}. Our use is consistent with~\cites{AP-tht,AHK}.

The above definition naturally extends to any polyhedral complex $Y$ and, for a face $\delta$ of $Y$, we denote by $\Sigma^\delta$ the corresponding star fan.

\medskip

We now recall what we mean by subdivisions and triangulations. We say that a polyhedral complex $Y$ is a subdivision of another polyhedral complex $Z$ if $Z$ and $Y$ share the same support, and in addition, each face of $Y$ is included in a face of $Z$. A subdivision $Y$ of $Z$ which is in addition simplicial is called a triangulation of $Z$.

\subsection{Canonical compactifications: local case} \label{sec:can-compact-fans-2}
We denote by $\eR = \R \cup \{\infty\}$ the extended real line with the topology extending that of $\R$ by a basis of open neighborhoods of $\infty$ given by intervals $(a, \infty]$ for any real number $a$. The addition of $\R$ naturally extends to $\eR$ and gives a monoid $(\eR, +)$ that we call the monoid of tropical numbers. We denote by $\eR_+ = \R_+ \cup\{\infty\}$ the submonoid of non-negative tropical numbers with the induced topology. Both monoids are modules over the semiring $\R_+$.

\medskip

For any cone $\sigma$ in $N_\R$, we denote by $\sigma^\vee \subseteq M_\R$ and $\sigma^\perp \subseteq M_\R$ the dual cone and the orthogonal plane to $\sigma$, respectively, defined by
\[\sigma^\vee := \Bigl\{m\in M_\R \st \langle m, a \rangle \geq 0 \ \textrm{\,for all\,}\ a \in \sigma\Bigr\}, \textrm{ and } \]
\[\sigma^\perp := \Bigl\{m \in M_\R \st \langle m, a \rangle = 0 \ \textrm{\,for all\,}\ a \in \sigma\Bigr\}. \]

\medskip

The canonical compactification $\comp \sigma$ of the cone $\sigma$, also called the extended cone associated to $\sigma$, is defined by the tensor product
\[\comp \sigma := \sigma \otimes_{\R_+} \eR_+, \]
endowed with the topology which is the finest one making all the endomorphisms
\[ z\mapsto z+z', \qquad a\mapsto x\otimes a, \qquad \textrm{and} \qquad x\mapsto x\otimes a \]
continuous for any $z,z'\in\comp\sigma$, $a\in\eR_+$ and $x\in\sigma$.
This topology turns the extended cone $\comp\sigma$ into a compact topological space. Moreover, its restriction to $\sigma$ coincides with the usual topology of $\sigma$.

We note that $\comp\sigma$ can be equivalently defined as follows. Denote by $\RpMod$ the category of $\R_+$-modules. Then, we have, as a set,
\[ \comp\sigma = \Hom_\RpMod(\sigma^\vee, \eR_+). \]

The extended cone $\comp\sigma$ has a distinguished point that we denote by $\infty_\sigma$. It is defined as $\infty_\sigma:=x\otimes\infty$ for any $x$ in the relative interior of $\sigma$. This definition does not depend on the choice of $x$. Alternatively, $\infty_\sigma$ is the point of $\Hom_\RpMod(\sigma^\vee, \eR_+)$ which takes value zero at any element of $\sigma^\perp$ and sends all the elements of $\sigma^\vee \setminus \sigma^\perp$ to $\infty$. Note in particular that we have $\infty_{\conezero} = 0$ where $\conezero$ denotes the cone $\{0\}$.

The definition of the canonical compactification is compatible with the inclusion of faces. For a face $\tau \subface \sigma$, we naturally get a map $\comp \tau \subseteq \comp \sigma$, which identifies $\comp \tau$ with the topological closure of $\tau$ in $\comp \sigma$.

\medskip

Let now $\Sigma$ be a rational fan in $N_\R$. We define the canonical compactification of $\Sigma$ denoted by $\comp \Sigma$ as the union of $\comp \sigma$ for any cone $\sigma$ in $\Sigma$, where the compactification $\comp \tau$ of $\tau$ is identified with the corresponding subspace of the extended cone $\comp \sigma$ for any $\tau \subface \sigma$ in $\Sigma$. We endow $\comp \Sigma$ with the induced quotient topology. In this way, each extended cone $\comp \sigma$ naturally embeds as a subspace of $\comp \Sigma$.

A rational fan $\Sigma$ gives naturally rise to a partial compactification of $N_\R$ that we denote by $\TP_\Sigma$. This partial compactification coincides with the tropicalization of the toric variety $\P_\Sigma$ associated to $\Sigma$. For this reason it is sometimes called the tropical toric variety defined by $\Sigma$.
The canonical compactification $\comp \Sigma$ naturally lives in $\TP_\Sigma$: in fact, it coincides with the closure of $\Sigma$ in $\TP_\Sigma$.

The tropical toric variety $\TP_\Sigma$ is defined as follows. We refer to~\cites{AP-hi, BGJK, Kaj08, OR11, Pay09, Thu07} for more details.

For any cone $\sigma$ in $\Sigma$, we set
\[ \~\sigma := \Hom_\RpMod(\sigma^\vee, \eR)\]
endowed with a natural topology that we do not precise here. Clearly, we have a natural inclusion of $\comp\sigma$ into $\~\sigma$. We set $N^\sigma_{\infty,\R}:=\infty_\sigma + N_\R \subseteq \~\sigma$. This is the set of all elements of $\Hom_\RpMod(\sigma^\vee, \eR)$ which take value $\infty$ at any point of $\sigma^\vee \setminus \sigma^\perp$ (and finite values on $\sigma^\perp$). Clearly $N^\conezero_{\infty,\R}=N_\R$. Moreover, the natural map
\[ N_\R\to N^\sigma_{\infty,\R},\quad z \mapsto z+\infty_\sigma \]
identifies $N^\sigma_{\infty,\R}\simeq N^\sigma_{\R}$. The space $\~\sigma$ admits a natural stratification into a disjoint union of subspaces $N^\tau_{\infty,\R} \simeq N^\tau_\R$, for $\tau$ a face of $\sigma$.

For a pair of cones $\tau \subface \sigma$ in $\Sigma$, we get an open inclusion $\~\tau\subseteq\~\sigma$. This allows to define $\TP_\Sigma$ as the gluing of $\~\sigma$, for $\sigma\in \Sigma$, along these inclusions. From the description above, we get a stratification of $\TP_\Sigma$ into the disjoint union of $N^\sigma_{\infty,\R} \simeq N^{\sigma}_\R$ for $\sigma \in \Sigma$.

For a point $x\in \TP_\sigma$ which lies in the stratum $N^\tau_{\infty,\R}$ for $\tau \subface \sigma$, the sedentarity of $x$ is by definition $\sed(x):=\tau$. In particular, if $\sigma$ is the positive quadrant in $\R^n$, the compactification $\TP_\sigma$ can be identified with $\eR^n$ and the faces of $\sigma$ with the subsets $[n]$ such that under these identifications, the sedentarity of a point becomes the subset of $[n]$ corresponding to those coordinates which are equal to $\infty$. This is consistent with the notation in~\cite{JSS19}.

\medskip

The canonical compactification of a fan $\Sigma$ admits a similar natural stratification into cones and fans that we describe now.

For a cone $\sigma \in \Sigma$ and a face $\tau \subface \sigma$ of $\sigma$, define $C^\tau_\sigma:=\infty_\tau+\sigma \subseteq \comp\sigma$. This coincides with the set of all elements of $\Hom_\RpMod(\sigma^\vee, \eR)$ which take value $\infty$ on $\tau^\vee \setminus \tau^\perp$, and finite values elsewhere. The cone $C^\tau_\sigma$ is isomorphic to the projection of $\sigma$ in $\rquot{N_{\sigma, \mathbb R}}{N_{\tau, \mathbb R}} \hookrightarrow N^{\tau}_{\infty,\R}$. We denote by $\Cint^\tau_\sigma$ the relative interior of $C^\tau_\sigma$.

For $\tau \in \Sigma$, the collection of cones $C^\tau_\sigma$ for $\sigma \in \Sigma$ with $\sigma \supface \tau$ form a fan in $N^\tau_{\infty, \R}$, with origin $\infty_\tau$, that we denote by $\Sigma_\infty^\tau\subseteq{N^\tau_{\infty,\R}}$. We note that the fan $\Sigma_\infty^\tau$ is canonically isomorphic to the star fan $\Sigma^\tau$. For the cone $\conezero$ of $\Sigma$, we have $\Sigma_\infty^\conezero = \Sigma$.

For any pair $\tau \subface \sigma$ in $\Sigma$, the closure $\comp C^\tau_\sigma$ of the cone $C^\tau_\sigma$ in $\comp \Sigma$ is the union of $\Cint^{\tau'}_{\sigma'}$ with $\tau\subface\tau'\subface\sigma'\subface\sigma$. Moreover, the closure of $\Sigma_\infty^\tau$ becomes canonically isomorphic to the canonical compactification of the fan $\Sigma^\tau \subseteq N^{\tau}_\R$.

\subsection{Extended polyhedral structures} Tropical compactifications of fans form examples of extended polyhedral complexes, which provide an enrichment of the category of polyhedral complexes and polyhedral spaces. We give a brief description and refer to~\cites{JSS19, MZ14, IKMZ, AP-geom} for more details. We restrict to the rational case since these are the only spaces we consider in this paper.

\medskip

Let $\sigma$ be a rational cone in $N_\R$. An extended polyhedron $\delta$ in $\TP_\sigma$ is the topological closure in $\TP_\sigma$ of any polyhedron included in a strata $N^\tau_{\infty,\R}$ for some $\tau\subface\sigma$. The topological closure of a face of $\delta \cap N^\zeta_{\infty,\R}$ for some face $\zeta$ of $\sigma$ is called a face of $\delta$. By an extended polyhedral complex in $\TP_\sigma$ we mean a finite collection $X$ of extended polyhedra in $\TP_\sigma$ verifying the two following properties:
\begin{itemize}
\item Any face $\gamma$ of an element $\delta\in X$ belongs to $X$.
\item The intersection of a pair of elements $\delta$ and $\delta'$ of $X$ is either empty or a common face of $\delta$ and $\delta'$.
\end{itemize}
The support of $X$ denoted by $\supp X$ is the union of $\delta \in X$. The space $\X = \supp X$ is called an extended polyhedral subspace of $\TP_\sigma$, and $X$ an extended polyhedral structure on $\X$.

\medskip

More general extended polyhedral spaces are then defined by using extended polyhedral subspaces of partial compactifications of vector spaces of the form $\TP_\sigma$ as local charts.

An integral extended polyhedral space $\X$ is by definition a Hausdorff topological space endowed with a finite atlas of charts $\Bigl(\phi_i\colon W_i\to U_i\subseteq\X_i\Bigr)_{i\in I}$, $I$ a finite set, with the following properties:
\begin{itemize}
\item The collection $\bigl\{\,W_i \st {i\in I}\, \bigr\}$ gives an open covering of $\X$.
\item Each $\X_i$, $i\in I$, is an extended polyhedral subspace of $\TP_{\sigma_i}$ for a finite dimensional real vector space $N_{i,\R}$ with $N_i$ a free $\Z$-module of finite rank and $\sigma_i$ a rational cone in $N_{i,\R}$. Moreover, $U_i$ is an open subset of $\X_i$.
\item The map $\phi_i$ is a homeomorphism between $W_i$ and $U_i$. Moreover, for any pair of indices $i,j\in I$, the transition map
\[\phi_j\circ \phi_i^{-1} \colon \phi_i(W_i \cap W_j) \to \TP_{\sigma_j}\]
is an extended integral affine map (with respect to lattices $N_i$ and $N_j$).
\end{itemize}
Recall that for two finite rank lattices $N_1$ and $N_2$, and rational cones $\sigma_1$ and $\sigma_2$ in $N_{1, \R}$ and $N_{2, \R}$, respectively, an extended integral affine map from an open subset $U \subseteq \TP_{\sigma_1}$ to $\TP_{\sigma_2}$ is a map which can be obtained as an extension of an integral affine map $\psi\colon N_{1, \R} \to N_{2, \R}$. By this we mean the following. Denote by $A$ the linear part of $\psi$ which is thus a $\Z$-linear map from $N_1$ to $N_2$. Consider the set $J$ consisting of all the rays $\varrho$ of $\sigma_1$ with the property that $A\varrho$ lives inside $\sigma_2$. Denote by $\tau_J$ the face of $\sigma_1$ generated by these rays. The affine map $\psi$ naturally extends to a map $\comp\psi\colon \bigcup_{\substack{\zeta\subface\tau_J}} (N_1)^\zeta_{\infty,\R} \to \TP_{\sigma_2}.$ We call this the extended affine map.

The requirement in the above definition is that there exists an extended integral affine map $\comp\psi$ from an open subset of $\TP_{\sigma_i}$ to $\TP_{\sigma_j}$ such that $\phi_i(W_i \cap W_j) \subseteq \bigcup_{\zeta\subface\tau_J} (N_1)^\zeta_{\infty,\R}$ for the corresponding face $\tau_J$ of $\sigma_i$, and the transition map $\phi_j\circ \phi_i^{-1}$ is the restriction of $\comp\psi$ to $\phi_i(W_i \cap W_j)$.

\medskip

A face structure on an extended polyhedral space $\X$ endowed with the corresponding atlas of charts, as above, is the choice, for each $i$, of an extended polyhedral complex structure $X_i$ with $\supp{X_i} = \X_i$, and a finite collection $\theta_1, \dots, \theta_N$ of closed set called facets for some integer $N \in \mathbb N$ such that the following properties hold:
\begin{itemize}
\item The facets cover $\X$.
\item Each $\theta_k$ is entirely contained in some chart $W_i$ for $i\in I$ so that the image $\phi_i(\theta_k)$ is the intersection of a face $\eta_{k,i}$ of $X_i$ with the open set $U_i$.
\item For a subset $J \subseteq [N]$ containing $k$, and for any chart $W_i$ containing $\theta_k$, the image of the intersection $\bigcap_{j\in J} \theta_j$ by $\phi_i$ in $U_i$ coincides with the intersection of a face of $\eta_{k,i}$ with $U_i$.
\end{itemize}

A face in this face structure is the preimage by $\phi_i$ of a face of $\eta_{j,i}$ for a $j\in[N]$ and for an $i\in I$ with $\theta_j\subseteq W_i$. Note that each face is contained in a chart $W_i$. We define the sedentarity of a face $\delta$ in a given chart $W_i$ as the sedentarity of any point in the relative interior of $\phi_i(\delta)$, viewed in $X_i$.

\subsection{Canonical compactifications: global case} \label{sec:can-compact-polycomp} Consider a polyhedral complex $Y$ in a real vector space $N_\R \simeq \R^n$. The recession pseudo-fan of $Y$ that we denote by $Y_\infty$ is the set of cones $\{\delta_\infty\st\delta\in Y\}$. In the case where this collection forms a fan, we call $Y_\infty$ the recession fan of $Y$. Any polyhedral complex $Y$ admits a subdivision whose recession pseudo-fan is a fan~\cite{AP-tht, OR11}.

Let now $Y$ be a polyhedral complex in $N_\R \simeq \R^n$ with recession fan $Y_\infty$. The canonical compactification of $Y$ denoted by $\comp Y$ is defined as the closure of $Y$ in the tropical toric variety $\TP_{Y_\sminfty}$.
It has a natural stratification given by cones $\sigma \in Y_\infty$. Consider a stratum $N^{\sigma}_{\infty, \R}$ of $\TP_\sigma$ and define $Y^\sigma_{\infty}$ as the intersection of $\comp Y$ with $N^{\sigma}_{\infty, \R}$. We drop $\infty$ if there is no risk of confusion and simply write $Y^\sigma$. We call $Y^{\conezero} = Y$ the open part of the compactification $\comp Y$. The boundary at infinity $D$ is defined as $D = \comp Y \setminus Y$. For each non-zero cone $\sigma$ in $Y_\infty$, let $D^\sigma$ be the closure of $Y^\sigma$ in $\comp Y$.

\begin{thm} [Tropical orbit-stratum correspondence]\label{thm:orbit-stratum-correspondence}
Notations as above, for each cone $\sigma \in Y_\infty$, we have the following.
\begin{enumerate}
\item The stratum $Y^\sigma$ is a polyhedral complex in $N^{\sigma}_{\infty, \R}$.
\item The recession pseudo-fan $(Y^{\sigma})_\infty$ of $Y^\sigma$ is a fan. Moreover, it coincides with the fan $(Y_\infty)_\infty^{\sigma}$ in $N_{\infty, \R}^\sigma$.
\item $D^\sigma$ coincides with the canonical compactification of $Y^{\sigma}$ in $N^{\sigma}_{\infty, \R}$, \ie, $D^\sigma = \comp{Y^{\sigma}}$.
\item If $Y$ has pure dimension $d$, then $Y^{\sigma}$ and $D^\sigma$ are of pure dimension $d -\dims \sigma$.
\end{enumerate}
\end{thm}
It follows from the above properties that $\comp Y$ is an extended polyhedral structure with a face structure induced from that of $Y$.

\subsection{Smooth tropical varieties} \label{sec:tropvar} Smoothness in tropical geometry reflects in polyhedral geometry the idea of maximal degeneracy for varieties defined over non-Archimedean fields. We will elaborate on this in our forthcoming work. Since this is a local notion, it amounts to fixing a good class of fans and their supports as local charts. Building on the Hodge theoretic interpretation of maximal degeneracy~\cite{Del-md}, our work~\cite{AP-tht} suggests these are fans which should satisfy the Poincar\'e duality for tropical cohomology, and which should have canonical compactifications of Tate type, meaning that the tropical cohomology of the canonical compactification is concentrated in Hodge bidegrees $(p,p)$ (see below for the definition of tropical cohomology). By~\cite{AP-hi}, the class of tropically smooth fans contains all Bergman fans including therefore complete fans. In this paper, Bergman fans are the ones which will serve as local charts, so we recall what we mean by a Bergman fan. 

A Bergman fan is a fan which has the same support as the Bergman fan of a matroid (the condition concerns only the support and not the fan structure itself).  Let $\Ma$ be a simple matroid on a ground set $E$ of rank $r+1$. We refer to~\cites{Oxl11, Wel10} for the definition and basic properties of matroids. A discussion can be also found in our paper~\cite{AP-tht}. Denote by $\{\e_i\}_{i\in E}$ the standard basis of $\Z^E$. For a subset $A \subseteq E$, let $\e_A$ be the sum $\sum_{i\in A} \e_i$ in $\Z^E$. Consider the lattice $N = \rquot{\Z^E}{\Z \e_E}$. The Bergman fan of $\Ma$ denoted by $\Sigma_\Ma$ is a rational fan in $N_\R$ of dimension $r$ defined as follows. A flag of proper flats $\Fl$ in $\Ma$ is a collection
\[\Fl\colon\quad\emptyset\neq F_1\subsetneq F_2\subsetneq \dots \subsetneq F_{\ell} \neq E\]
consisting of proper flats $F_1, \dots, F_\ell$ of $\Ma$.
To such a flag, we associate the rational cone $\sigma_\Fl \subseteq N_\R$ of dimension $\ell$ generated by the vectors $\e_{F_1}, \e_{F_2}, \dots, \e_{F_\ell}$, \ie,
\[\sigma_\Fl := \Bigl\{\lambda_1\e_{F_1} + \dots + \lambda_\ell \e_{F_\ell} \Bigst \lambda_1, \dots, \lambda_\ell \geq0\Bigr\}.\]
The Bergman fan $\Sigma_\Ma$ of $\Ma$ consists of all the cones $\sigma_\Fl$ for $\Fl$ a flag of proper flats of $\Ma$:
\[\Sigma_\Ma := \Bigl\{\, \sigma_\Fl \Bigst \Fl \textrm{ flag of proper flats of }\Ma\,\Bigr\}.\]
It has pure dimension $r$.

A Bergman support in a real vector space $W$ is a subset $S\subseteq W$ which is isomorphic to $\supp{\Sigma_\Ma}$ via a linear map $\phi\colon W \to \rquot{\R^E}{\R \e_E}$. A fan $\Sigma$ in $W$ is called Bergman if its support $\supp{\Sigma}$ is Bergman. If $W$ is equipped with a full rank lattice $N_W$ and the map $\phi$ induces an isomorphism between the two lattices $N_W\cap \TT S$ and $N$, we say $\Sigma$ is a rational Bergman fan.

Note that the above terminologies are consistent in the sense that the Bergman fan $\Sigma_\Ma$ of a matroid $\Ma$ is an example of a Bergman fan. In addition, any complete fan in a real vector space is Bergman. The tropicalization of the complement of a hyperplane arrangement is Bergman as well~\cite{AK06}. The category of Bergman fans is closed under product~\cite{AP-hi}.

\smallskip
A smooth tropical variety is an extended polyhedral space with an integral affine structure that is locally modeled by supports of Bergman fans~\cites{IKMZ, MZ14, JSS19}. In other words, any point has a neighborhood which is isomorphic to an open set in $\Sigma\times\eR^k$ for some Bergman fan $\Sigma$.

We have the following theorem~\cite{AP-geom}.

\begin{thm}\label{thm:smoothness-compact} Consider a rational polyhedral complex $Y$ in $N_\R$ with smooth support. Assume the recession fan $Y_\infty$ of $Y$ is unimodular. In this case, the canonical compactification $\comp Y$ of $Y$ has smooth support. If the polyhedral structure on $Y$ is unimodular, then the extended polyhedral structure induced on $\comp Y$ is unimodular as well.
\end{thm}
In particular, canonical compactifications of unimodular Bergman fans are smooth.

\medskip

A rationally triangulable smooth projective tropical variety $X$ is by definition a smooth tropical variety which is isomorphic to the canonical compactification $\comp Y$ of a rational polyhedral complex $Y$ with smooth support such that the recession fan $Y_\infty$ is unimodular and quasi-projective. Recall that a rational fan $\Sigma$ in $N_\R$ is called quasi-projective if the corresponding toric variety $\P_\Sigma$ is quasi-projective~\cite{Ful-toric}. By the triangulation theorem proved in~\cite{AP-tht}, any rationally triangulable smooth projective tropical variety admits a triangulation which is unimodular with respect to the lattice $\frac 1m N$ for some positive integer $m$.

\subsection{Tropical homology and cohomology} Tropical homology and cohomology groups were introduced by Itenberg-Katzarkov-Mikhalkin-Zharkov~\cites{IKMZ} and further studied in~\cites{JSS19, MZ14, JRS19, GS-sheaf, AP-tht}. We briefly recall the definition of these groups.

Consider an extended polyhedral space $X$ with a face structure. We define the multi-tangent and multi-cotangent coefficient groups $\SF_p(\delta)$ and $\SF^p(\delta)$ associated to each face $\delta$ of $X$. All together, they lead to the definition of chain and cochain complexes which define the tropical homology and cohomology groups of $X$.

Let $\delta$ be a face of $X$ and $p$ be a non-negative integer. The $p$-th multi-tangent and the $p$-th multi-cotangent spaces $\SF_p(\delta)$ and $\SF^p(\delta)$ of $X$ at $\delta$ are defined by
\[\SF_{p}(\delta)= \hspace{-.5cm} \sum_{\eta \supface \delta \\\sed(\eta)=\sed(\delta) } \hspace{-.5cm} \bigwedge^p\TT\eta,\qquad \textrm{and} \qquad \SF^p(\delta) = \SF_p(\delta)^\dual, \]
where as before $\SF_p(\delta)^\dual$ means the dual of $\SF_p(\delta)$. An inclusion of faces $\gamma \subface \delta$ of $X$ gives natural maps $\i_{\delta\supface\gamma}\colon \SF_p(\delta)\to\SF_p(\gamma)$ and $\i^*_{\gamma\subface\delta}\colon \SF^p(\gamma) \to \SF^p(\delta)$.

\medskip

For a non-negative integer $p$, the cellular chain complex
\[C_{p,\bul}\colon \quad \dots\longrightarrow C_{p, q+1}(X) \xrightarrow{\partial^\trop_{q+1}} C_{p,q}(X) \xrightarrow{\ \partial^\trop_{q}\ } C_{p,q-1} (X)\longrightarrow\cdots\]
is defined by setting
\[C_{p,q}(X) := \bigoplus_{\delta \in X \\ \dims{\delta} =q} \SF_p(\delta), \]
for any non-negative integer $q$, and by using maps $\i_{\delta\supface\gamma}$ with signs as in cellular homology theory. The tropical homology of $X$ is defined by taking the homology of the tropical chain complex, that is,
\[H_{p,q}^\trop(X) := H_q(C_{p,\bul}).\]
Similarly, we can define the cochain complex
\[C^{p,\bul}\colon \quad \dots\longrightarrow C^{p, q-1}(X) \xrightarrow{\partial_\trop^{q-1}} C^{p,q}(X) \xrightarrow{\ \partial_\trop^{q}\ } C^{p,q+1}(X) \longrightarrow\cdots\]
with
\[C^{p,q}(X) := C_{p,q}(X)^\dual \simeq \bigoplus_{\delta\in X \\
\dims{\delta}=q} \SF^p(\delta),\]
and define the tropical cohomology of $X$ by
\[H^{p,q}_\trop(X) := H^q(C^{p,\bul}).\]

If the polyhedral structure $X$ is rational, tropical homology and cohomology groups can be defined with integer or rational coefficients. In fact, for each face $\delta$, we have $\TT\delta = N_{\delta, \R}$ for the lattice $N_\delta$ associated to $\delta$, and one can define
\[\SF_{p}(\delta,\Z)= \!\!\sum_{\eta \supface \delta \\\sed(\eta)=\sed(\delta)}\!\!\bigwedge^p N_\eta, \qquad \textrm{and} \qquad \SF^{p}(\delta, \Z) = \SF_p(\delta,\Z)^{\vee}\]
and define the corresponding complexes $C_{p,\bul}^\Z$ and $C^{p, \bul}_\Z$ with $\Z$-coefficients. This gives
\[H_{p,q}^\trop(X, \Z) := H_q(C_{p,\bul}^{\Z}) \qquad H^{p,q}_\trop(X, \Z) := H^q(C^{p,\bul}_\Z).\]
Similarly, we get $H_{p,q}^\trop(X, \Q)$ and $H^{p,q}_\trop(X, \Q)$.

\smallskip
Compact smooth tropical varieties satisfy Poincar\'e duality. This was proved by Jell-Shaw-Smacka~\cite{JSS19} with rational coefficients and by Jell-Rau-Shaw~\cite{JRS19} for integral coefficients.

\medskip

In this paper, unless otherwise stated, the cohomology and homology groups are all with rational coefficients.

\section{Minkowski weights and tropical cycles}\label{sec:minkowski_weight}

The aim of this section is to introduce the Minkowski weights and to explain how they give rise to tropical cycles. We refer to~\cites{AR10, MS15, MR09, GS19} for the definition of tropical cycles and their basic properties.

\subsection{Minkowski weights} Let $p, d$ be two non-negative integers with $p\leq d$, and let $Y$ be a unimodular polyhedral complex with smooth support of pure dimension $d$ in $N_\R$ with a unimodular recession fan $Y_\infty$. By Theorem~\ref{thm:smoothness-compact}, the closure $X=\comp Y$ of $Y$ in $\TP_{Y_{\infty}}$ is smooth. We denote by $\X$ and $\Y$ the supports of $X$ and $Y$, respectively.

Assume for each face $\delta$ of $Y$ of dimension $d-p$ we are given a weight which is an integer (or a rational number, depending on the context) denoted by $w(\delta)$. Let $C := (Y_{(d-p)}, w)$ be the corresponding weighted polyhedral complex $Y$ with the weight function $w$ on the facets of $Y_{(d-p)}$. The weight function $w$ is called a \emph{Minkowski weight} of dimension $d-p$ on $Y$ if the following \emph{balancing condition} is verified:
\[\forall\:\gamma\in Y_{d-p-1},\qquad \sum_{\delta\ssupface\gamma} w(\delta) \e_{\delta/\gamma} = 0 \in N^{\gamma}.\]
Here $\e_{\delta/\gamma}$ is the primitive vector of the ray $\rho_{\delta/\gamma}$ corresponding to $\delta$ in $N^\gamma_\R$ and the sum is over the face $\delta$ such that $\gamma$ is a subface of codimension one in $\delta$. We denote by $\MW_{d-p}(Y)$ the set of all Minkowski weights of dimension $d-p$ on $Y$. Addition of weights cell by cell turns $\MW_{d-p}(Y)$ into a group.

From an element of $w\in\MW_{d-p}(Y)$, we naturally obtain an element of the homology group $H_{d-p,d-p}^\trop(X)$ by setting on each $(d-p)$-dimensional face $\delta$ the canonical unit multivector $\nvect_\delta$ of $\delta$ in $F_{d-p}(\delta)$ with the corresponding coefficient $w(\delta)$.

\subsection{Tropical cycle associated to a Minkowski weight} Each Minkowski weight on $Y$ gives a tropical cycle in $\X$. This is obtained by taking the closure $\comp C = (\comp Y_{(d-p)}, w)$ of $Y_{d-p}$ in $\X$ with the same weight function $w$ on its facets. We call $\comp C$ the tropical cycle in $X$ corresponding to the Minkowski weight $w$ on $Y$. Tropical cycles in $\X$ which are of this form, for a given choice of a rational polyhedral structure $Y$ on the open part $\Y$ of $\X$, are called \emph{admissible}.

\subsection{Relation with tropical homology} Each tropical cycle of codimension $p$ gives by integration an element of $H_\Dolb^{d-p,d-p}(X, \R)^\dual$, and thus of $H_\trop^{d-p,d-p}(X, \R)^\dual$ via the isomorphism between tropical Dolbeault cohomology and tropical singular cohomology \cite{JSS19}. This is depicted in the diagram of Figure \ref{fig:global_diagram}. In this diagram, we represent the set of tropical cycles of codimension $p$ on $\X$ by $\Cycle_{d-p}(\X)$.
\begin{figure}[t]
\[ \begin{tikzcd}[column sep=small, row sep=small]
\MW_{d-p}(Y) \rar\dar &
H^\trop_{d-p,d-p}(X) \rar &
H^{d-p,d-p}_\trop(X)^\dual \rar["\simeq", "\PD"'] \dar[hook] &
H^{p,p}_\trop(X)
\\
\Cycle_{d-p}(\X) \rar &
H^{d-p,d-p}_\Dolb(\X, \R)^\dual \rar{\simeq} &
H^{d-p,d-p}_\trop(X, \R)^\dual
\end{tikzcd} \]
\caption{Main relations between Minkowski weights, tropical cycles and tropical homology and cohomology. \label{fig:global_diagram}}
\end{figure}

\section{Integral tropical Hodge conjecture in the local case} \label{sec:local_HC}

In this section, unless otherwise explicitly stated, we work with integral coefficients.

\subsection{Hodge isomorphism for Bergman fans} Let $\Sigma$ be a unimodular Bergman fan of dimension $d$ and consider its canonical compactification $\comp\Sigma$. In this section we explain why Theorem \ref{thm:HI} induces an integral version of the Hodge conjecture.

\begin{thm}[Hodge isomorphism for Bergman fans] \label{thm:HI} For a unimodular Bergman fan $\Sigma$ of dimension $d$, the cycle class map
\[\class\colon A^p(\Sigma) \longsimto H^{p,p}_\trop(\comp \Sigma)\]
induces an isomorphism between the Chow groups of $\Sigma$ and the tropical cohomology groups of $\comp\Sigma$. Moreover, this isomorphism is compatible with the diagram of Figure \ref{fig:local_diagram}.
\end{thm}

This theorem can be regarded as the tropical analogue of a theorem of Feichtner and Yuzvinsky~\cite{FY04}, which establishes a similar result for wonderful compactifications of hyperplane arrangement complements. We furthermore show in~\cite{AP-hi} that the cohomology groups $H_\trop^{p,q}(\comp \Sigma)$ for $p\neq q$ are all trivial.

\subsection{The Chow ring of a Bergman fan} \label{sec:chow_ring} We recall the definition of the Chow groups; connection to tropical cycles is explained in the next subsection. For any ray $\varrho\in\Sigma_1$, let $\e_\varrho$ be the generator of $\varrho\cap N$.

The Chow ring $A^\bul(\Sigma)$ is the graded ring defined as the quotient of the polynomial ring $\Z[\x_\varrho \mid \varrho\in \Sigma_1]$, with generator $\x_{\varrho}$ associated to the ray $\varrho$ of $\Sigma$, by the homogeneous ideal $\II_1 +\II_2$ where \begin{itemize}
\item the ideal $\II_1$ is generated by monomials of the form $\prod_{\varrho \in S} \x_\varrho$ for any subset $S \subseteq \Sigma_1$ of rays which does not form a cone in $\Sigma$; and
\item the ideal $\II_2$ is generated by elements of the form
\[\sum_{\varrho\in \Sigma_1} \langle m, \e_{\varrho}\rangle \x_\varrho\]
for any element $m$ in the dual lattice $M = N^\vee$.
\end{itemize}
Here $\langle \ccdot, \rdot\rangle$ is the duality pairing between $M$ and $N$.

\medskip

If $\sigma$ is any cone of $\Sigma$ of dimension $p$, we set
\[ \x_\sigma := \prod_{\varrho \in \Sigma_1 \\ \varrho \subface \sigma} \x_\varrho \in A^p(\sigma). \]
One can show that for any two maximal cones $\eta, \eta' \in \Sigma_d$, we have $\x_\eta = \x_{\eta'}$. Moreover, the top degree part $A^d(\Sigma)$ is of rank one and is generated by $\x_\eta$.

From this, we deduce the degree map $\deg\colon A^d(\Sigma) \simto \Z$ by mapping $\alpha$ onto $1$. This map induces a perfect pairing
\[ \arraycolsep=1.4pt \begin{array}{rclcl}
A^k(\Sigma) & \times & A^{d-k}(\Sigma) & \to & \Z, \\
y & , & z & \mapsto & \deg(y \cdot z).
\end{array} \]

\subsection{Integral tropical Hodge conjecture in the local case} We can adapt the diagram of Figure \ref{fig:global_diagram} to our local case. The diagram can be furthermore completed thanks to the following theorem from~\cite{AHK}.

\begin{thm}\label{thm:hc-local} Let $\Sigma$ be a unimodular fan in $N_\R$ of dimension $d$. Then, there is an isomorphism $A^p(\Sigma) \simeq \MW_{d-p}(\Sigma)$.
\end{thm}

Thereby, we obtain the diagram in Figure \ref{fig:local_diagram}. By Theorem \ref{thm:HI}, every map in the first row is an isomorphism. (The one in the middle concerning $H^\trop_{d-p,d-p}(\comp\Sigma, \Z) \simeq
H^{d-p,d-p}_\trop(\comp\Sigma, \Z)^\vee$ is explained in~\cite{AP-hi}.) As a consequence, to any element $\alpha$ of $H^{p,p}_\trop(\comp\Sigma, \Z)$ we can associate an admissible tropical cycle in $\Cycle_{d-p}(\comp\Sigma)$ with integral coefficients whose image by the cycle class map is $\alpha$. This is the integral Hodge conjecture in the local case.

\begin{figure}[t]
\[ \begin{tikzcd}[column sep=small, row sep=small]
A^p(\Sigma, \Z) \rar{\simeq} &
\MW_{d-p}(\Sigma, \Z) \rar\dar &
H^\trop_{d-p,d-p}(\comp\Sigma, \Z) \rar &
H^{d-p,d-p}_\trop(\comp\Sigma, \Z)^\vee \rar["\simeq", "\PD"'] \dar[hook] &
H^{p,p}_\trop(\comp\Sigma, \Z)
\\
&
\Cycle_{d-p}(\comp\Sigma, \Z) \rar &
H^{d-p,d-p}_\Dolb(\comp\Sigma, \R)^\dual \rar{\simeq} &
H^{d-p,d-p}_\trop(\comp\Sigma, \R)^\dual
\end{tikzcd} \]
\caption{The local analogue of the diagram of Figure \ref{fig:global_diagram}. \label{fig:local_diagram}}
\end{figure}

\medskip

We will prove Theorem~\ref{thm:HC} by generalizing this picture to the global setting, by showing that each element in the kernel of the monodromy is the image by the cycle class map of an admissible tropical cycle with rational coefficients.

\section{Tropical Steenbrink double complex and the comparison theorem} \label{sec:steenbrink}

Let $\X$ be a rationally triangulable compact smooth tropical variety of dimension $d$. Let $X$ be a unimodular triangulation of $\X$. Denote by $X_\f$ the set of faces of $X$ that do not intersect the boundary at infinity of $X$. In~\cite{AP-tht} we defined the first page of the tropical Steenbrink sequence associated to $X$ by
\[ \ST_1^{a,b} := \bigoplus_{s\geq\abs a \\ s\equiv a\pmod 2} \ST_1^{a,b,s} \]
where
\[ \ST_1^{a,b,s} := \bigoplus_{\delta\in X_\f \\ \dims\delta = s}H^{a+b-s}(\delta). \]
Here, we set $H^{k}(\delta):=H^{k}(\comp \Sigma^\delta) = \bigoplus_{p+q=k} H^{p,q}_\trop(\comp\Sigma^\delta)$. The definition is motivated by classical Hodge theory where the Steenbrink spectral sequence is a spectral sequence which calculates the weight-graded pieces of the limit mixed Hodge structure for a degenerating family $\fX^*$ over the punctured disk $\disk^*$. It has a first page with a shape similar to the tropical one above, and described by the the special fiber of a semistable extension $\fX$ of the family $\fX^*$ over $\disk$ (which exists after a finite ramified base change of the base). In the tropical setting, the group $H^{k}(\delta)$ is zero unless $k$ is even in which case, we get $H^{k}(\delta) = H^{k/2,k/2}_\trop(\comp\Sigma^\delta) \simeq A^{k/2}(\Sigma^\delta)$. In particular, the bigraded piece $\ST_1^{a,b}$ in the tropical setting is trivial if $b$ is odd.

In the tropical Steenbrink sequence, the differentials of bidegree $(1,0)$ are given by $\d:=\gys+\i^*$, where
\[\i^{a,b\,*} \colon \ST^{a,b}_1 \to \ST_1^{a+1, b}, \qquad \textrm{and} \qquad \gys^{a,b}\colon \ST^{a,b}_1 \to \ST_1^{a+1, b}.\]
The map $\i^*$ corresponds to the restriction in cohomology, and $\gys$ is the Poincar\'e dual of the restriction map. The precise definition can be found in~\cite{AP-tht}. We proved in~\cite{AP-tht} that for a unimodular triangulation $X$ of $\X$ and for any integer $b$, the differential $\d$ makes $\ST_1^{\bul,b}$ into a cochain complex. We denote by $H^{a}(\ST^{\bul,b}_1, \d)$ the $a$-th cohomology of this cochain complex.

\smallskip
The cohomology of the Steenbrink cochain complex can be computed by the following theorem.
\begin{thm}[Steenbrink-Tropical comparison theorem] \label{thm:steenbrink-tropical} Notations as above, let $\X$ be a compact smooth tropical variety of dimension $d$. Let $X$ be a unimodular triangulation of $\X$. The cohomology of $(\ST_1^{\bul,b},\d)$ is described in the following way. If $b$ is odd, then all the terms $\ST_1^{a,b}$ are zero, and the cohomology is vanishing. For $b$ even, writing $b=2p$ for $p \in \mathbb Z$, we have for any $q\in \mathbb Z$ a canonical isomorphism
\[H^{q-p}(\ST^{\bul,2p}_1, \d) \simeq H^{p,q}_{\trop}(\X).\]
\end{thm}

In the complex approximable case, i.e., when $\X$ arises as the tropicalization of a family of complex projective varieties, this theorem is proved in~\cite{IKMZ}. The statement in this generality is proved in~\cite{AP-tht}. Related work on integral affine manifold with singularities can be found in~\cites{GS10, GS06, Rud10, Rud20}.

\medskip

This isomorphism restricted to the kernel of the monodromy in $H^{p,p}_\trop(\X)$ will be explicitly described in Section \ref{sec:explicit_isomorphism}.

\subsection{Monodromy} The monodromy operator $N$ on $\ST_1^{\bul,\bul}$ is of bidegree $(2,-2)$ and is given by the identity map $\id$ on the relevant parts. More precisely, $N\colon \ST_1^{a,b}\to\ST_1^{a+2,b-2}$ is given on each nontrivial part by
\[ N\rest{\ST_1^{a,b,s}} = \begin{cases}
  \id\colon \ST_1^{a,b,s}\to \ST_1^{a+2,b-2,s} & \text{if $s\geq\abs{a+2}$}, \\
  0 & \text{otherwise.}
\end{cases} \]

It is not hard to see that it is a Lefschetz operator which verifies the Hard Lefschetz property around $0$, that is, the map
\[ N^k\colon \ST_1^{-k,b+k} \to \ST_1^{k,b-k} \]
is an isomorphism for any $k$. The Hard Lefschetz property around $0$ will be denoted by $\HL$ in the sequel.

The operator $N$ commutes moreover with the differential and induces a monodromy operator on the cohomology $\ST_2^{\bul,\bul}$ that we still denote by $N$. This induced operator verifies as well $\HL$ as we proved in \cite{AP-tht}.

\medskip

We have the following characterization of the monodromy operator in terms of the eigenwave map $\phi$ defined in~\cite{MZ14} and the monodromy on tropical Dolbeault cohomology~\cites{Liu19}.

\begin{thm} \label{thm:eigenwave_monodromy} The monodromy operator $N\colon H^{p,q}_\trop(\X) \to H^{p-1, q+1}_\trop(X)$ coincides with the eigenwave operator $\phi$. With real coefficients, it coincides as well with the monodromy operator on Dolbeault cohomology.
\end{thm}
\begin{proof} The proof of the first statement is given in~\cite{MZ14} in the complex approximable case, \ie, in the case the tropical variety arises as a tropical limit of a family of complex projective varieties, and it is proved in~\cite{AP-tht} for the general case. The proof of the second statement is a consequence of \cite{Jel19} which relates the eigenwave operator to the monodromy operator on tropical Dolbeault cohomology.
\end{proof}

\subsection{Polarization} We now describe a natural polarization on $\ST_1^{\bul, \bul}$, which after passing to the cohomology, induces a polarization on cohomology groups $H^{a}(\ST_1^{\bul, b}, \d)$. We refer to~\cite{AP-tht} for more details. The material in this section will allow us later to describe the cohomology classes associated to Minkowski weights.

\medskip

We can define a natural bilinear form $\psi$ on $\ST_1^{\bul, \bul}$ as follows. Take two elements $x \in \ST_1^{a,b,s}$ and $y \in \ST_1^{a',b',s'}$, for integers $a,b,s, a',b',s'$. We can write
\[ x=\sum_{\delta\in X_\f \\ \dims\delta=s} x_\delta \]
with $x_\delta\in H^{a+b-s}(\delta)$, and similarly for $y$. The degree map is defined by setting
\[ \deg(x \cdot y) := \begin{cases}
\sum_{\delta\in X_\f \\ \dims\delta=s}\deg(x_\delta\cdot y_\delta) & \text{if $s=s'$}, \\
0 & \text{otherwise}.
\end{cases} \]
where $\deg(x_\delta \cdot y_\delta)$ is defined in Section \ref{sec:chow_ring}. Note that the notation $x \cdot y$ in the above expression simply means $\sum_\delta x_\delta y_\delta$ whenever this makes sense, with the sum running over the corresponding $\delta$, and it is set to zero otherwise. The bilinear form $\psi$ is now defined by
\[ \psi(x,y) := \begin{cases}
\epsilon(a,b) \deg(x \cdot y) & \textrm{if $a+a'=0$, \ $b+b'=2d$, and $s'=s$,} \\
0 & \textrm{otherwise},
\end{cases} \]
where for a pair of integers $a,b$, with $b$ even, we set
\[\epsilon(a,b) := (-1)^{ a+ \frac b2},\]
and for $b$ odd, we define $\epsilon(a,b) :=1$. In this latter case, we necessarily have $x=0$.

\medskip

We proved in~\cite{AP-tht} that the bilinear pairing $\psi$ satisfies the following nice properties.

\begin{prop} \label{prop:properties_psi}
For any pair of elements $x, y \in \bigoplus_{a,b} \ST_1^{a,b}$, we have
\begin{enumerate}[label=\defaultRoman]
\item \label{enum-psi:i} $\psi(x,y) = (-1)^d \psi(y,x)$.
\item \label{enum-psi:ii} $\psi(Nx, y) + \psi(x, Ny) =0$.
\item \label{enum-psi:} $\psi(\d x, y) + \psi(x, \d y) =0$.
\end{enumerate}
\end{prop}

It follows from this proposition that the polarization $\psi$ induces a polarization on the cohomology groups $H^{a}(\ST_1^{\bul, b}, \d)$. By this we mean the following. First, since $[\d,N]=0$, we get an induced map
\[N\colon H^{a}(\ST_1^{\bul, b}, \d) \to H^{a+2}(\ST_1^{\bul, b-2}, \d).\]
Now, since $\psi(\d\ccdot,\rdot)=-\psi(\ccdot,\d\ccdot)$, we get an induced pairing
\[ \psi\colon H^{\bul}(\ST_1^{\bul, \bul}, \d)\times H^{\bul}(\ST_1^{\bul, \bul}, \d) \to \Q. \]
We have the following theorem.

\begin{thm} \label{thm:tropical_differential_HL_structure} The following properties hold. For any integer $a\geq 0$, we have
\begin{itemize}
\item The map $N^a \colon H^{-a}(\ST_1^{\bul, b}, \d) \to H^a(\ST_1^{\bul, b-2a}, \d)$ is an isomorphism. In other words, the induced monodromy operator $N$ on cohomology groups satisfies the Hard Lefschetz property $\HL$.
\item The polarization $\psi$ induces a non-degenerate bilinear form $\psi(\ccdot,N^a \ccdot)$ on $H^{-a}(\ST_1^{\bul, b}, \d) \times H^{-a}(\ST_1^{\bul, 2d-b+2a}, \d)$.
\item Let $a\geq 0$ and $b$ be two integers, and denote by $P^{-a,b}$ the primitive part of $H^{-a}(\ST_1^{\bul, b}, \d)$ with respect to the monodromy defined by
\[ P^{-a,b} := \ker\Bigl(N^{a+1} \colon H^{-a}(\ST_1^{\bul, b}, \d)\longrightarrow H^{a+2}(\ST_1^{\bul, b-2a-2}, \d)\Bigr).\]
Then, we have the decomposition into primitive parts
\[ H^{-a}(\ST_1^{\bul, b},\d) = \bigoplus_{s \geq 0} N^s P^{-a-2s,b+2s}.\]
Moreover, this decomposition is orthogonal for the pairing $\psi(\ccdot,N^{a}\ccdot)$ meaning that the term $N^s P^{-a-2s,b+2s}$ of the above decomposition and the term $N^{s'} P^{-a-2s',2d-b+2a+2s'}$ in the primitive decomposition of $H^{-a}(\ST_1^{\bul, 2d-b+2a}, \d)$ are orthogonal for $s\neq s'$.
\end{itemize}
\end{thm}

\begin{proof} This follows from the assumption that $\X$ is projective, which shows the existence of a tropical K\"ahler form on $\X$, and from Theorem 6.19 in~\cite{AP-tht} which proves a more refined statement in terms of the corresponding Hodge-Lefschetz structure.
\end{proof}

By Steenbrink-tropical comparison theorem, the above statement can be translated to tropical cohomology. In particular, we get the decomposition
\[H_\trop^{p,q}(\X, \Q) = \bigoplus_{s \geq 0} N^s PH_\trop^{p+s, q-s}(\X, \Q)\]
where for any pair of integers $p,q$ with $p\geq q$, we set
\[PH_\trop^{p, q}(\X, \Q) := \ker\Bigl(N^{p-q+1} \colon H^{p,q}_\trop(\X, \Q)\longrightarrow H_\trop^{q-1, p+1}(\X, \Q)\Bigr),\]
for the primitive parts with respect to the monodromy operator.
Moreover, the decomposition is orthogonal with respect to the pairing $\psi(\ccdot,N^{p-q}\ccdot)$.

\subsection{Explicit isomorphism between primitive parts $P^{0,2p}$ and $PH^{p,p}_\trop$} \label{sec:explicit_isomorphism} In this section we describe more explicitly  the isomorphism given by the comparison Theorem \ref{thm:steenbrink-tropical} restricted to primitive parts of respective bidegree $(0,2p)$ and $(p,p)$ in the Steenbrink sequence and tropical cohomology groups, respectively. These are the relevant parts for our purpose.

\subsubsection{Idea of the proof of the comparison theorem} Theorem \ref{thm:steenbrink-tropical} is proved in Section 5 of \cite{AP-tht}. We recall the main ideas here.

We wish to link the cohomology groups of the tropical Steenbrink sequence $\ST^{\bul,2p}_1$ with those of $C^{p,\bul}_\trop(X)$. The basic idea is based on the use of the following lemma.

\begin{lemma}[Zigzag isomorphism] \label{lem:zigzag}
Let $\AA^{\bul,\bul}$ be a double complex with differentials $\d$ and $\d'$ of respective degree $(1,0)$ and $(0,1)$. Assume that
\begin{itemize}
\item $\d\d'+\d'\d=0$,
\item $\AA^{b,b'}=0$ if $b<0$ or $b'<0$,
\item $\AA^{b,\bul}$ is exact if $b>0$,
\item $\AA^{\bul,b'}$ is exact if $b'>0$.
\end{itemize}
Then, there is a canonical isomorphism
\[ H^\bul(\AA^{0,\bul})\simeq H^\bul(\AA^{\bul,0}). \]
\end{lemma}

In this lemma, the isomorphism is given by $\d'^{-1}\d\d'^{-1}\cdots\d\d'^{-1}\d\colon \AA^{0,p} \to \AA^{p,0}$, where, by $\d'^{-1}$ we mean choosing any preimage by $\d'$ of an element. This map is not well-defined a priori. Nevertheless, we can prove that it induces a well-defined map on the level of the cohomology.

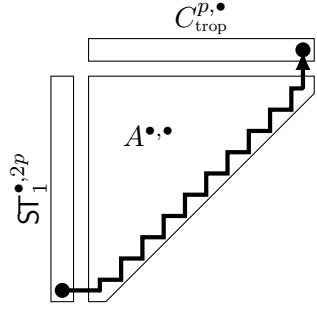
\begin{figure}[tp]
\begin{tikzpicture}[line cap=round]
\draw (0,.2) rectangle ++ (3,.3) node[midway, above=2pt] {$C^{p,\bul}_\trop$};
\draw (-.2,0) rectangle ++ (-.3,-3) node[midway, left=2pt] {\rotatebox{90}{$\ST_1^{\bul,2p}$}};
\draw (0,0) -- (3,0) --++ (0,-.23) -- (.23,-3) -- (0,-3) -- cycle;
\draw (.8,-.8) node {$A^{\bul,\bul}$};
\fill (-.35,-3) ++ (0,.15) coordinate (A) circle (.1);
\fill (3,.35) ++ (-.15,0) coordinate (B) circle (.1);
{
\draw[ultra thick] (A) -- (.15,-3+.15) coordinate (A');
\draw[ultra thick, latex-] (B) -- (3-.15,-.15) coordinate (B');
\draw[decorate, decoration={zigzag,segment length=4mm,amplitude=1mm}, ultra thick] (A') -- (B');}
\end{tikzpicture}
\caption{Illustration of the zigzag lemma between $\ST_1^{\bul,2p}$ and $C_\trop^{p,\bul}$ \label{fig:zigzag-scheme}}
\end{figure}

\smallskip
Coming back to our situation, assume now that there exists a double complex $\AA$ which verifies the conditions of the lemma and moreover, $\AA^{0,\bul}=C^{p,\bul}[-1]$ and $\AA^{\bul, 0}=\ST^{\bul,2p}_1[-1]$ as in Figure \ref{fig:zigzag-scheme}. Then the zigzag lemma implies the isomorphism in cohomology.

\smallskip
Actually, this is not exactly what we do in \cite{AP-tht} as the situation is more complicated. We first need to introduce spectral sequences that compute the tropical cohomology groups and the cohomology groups of $\ST^{\bul,2p}_1$. Then we are able to introduce a triple complex, which is the analogue in the setting of spectral sequences of the double complex $\AA^{\bul,\bul}$ above. This we can do thanks to the tropical Deligne resolution of the coefficient sheaves (cf. Section 5.3 of \cite{AP-tht}). Finally, we need to prove a more general version of the zigzag lemma, called the spectral resolution lemma, in order to link both these spectral sequences. All this turns out to be quite technical.

\smallskip
This being said, as far as for the parts which are of interest to us in this paper, everything happens to be luckily on the border, and the restriction of the isomorphism and the analysis can be made completely explicit. The idea is depicted in Figure \ref{fig:zigzag-scheme} and the details will be given in the next section. In that figure, the zigzag arrow is the part we have to study in (the spectral resolution analogue of) $\AA^{\bul,\bul}$; everything on the right hand side and below this zigzag arrow is trivial in the double complex. The arrow itself is given in detail in Figure \ref{fig:zigzag}. The start of the arrow is the kernel of the monodromy in $\ST_1^{0,2p}$. The end of the arrow is a subquotient of $C^{p,p}_\trop(X)$.

\smallskip

With this preparation, we can now explain the details of this scheme.

\subsubsection{Details of the zigzag} Looking at the definitions, we get that
\[ \ker(N\colon \ST_1^{0,2p} \to \ST_1^{2,2p-2}) = \bigoplus_{v \in X_{\f,0}} H^{2p}(v), \]
where the sum is over vertices of $X_\f$. Hence we need to describe the map from $\bigoplus_{v \in X_{\f,0}} H^{2p}(v)$ to $C^{p,p}_\trop(X)$ which induces the isomorphism on the cohomology given by the Steenbrink-Tropical comparison theorem.

\begin{figure}[th]
\renewcommand{\top}[2]{\displaystyle\bigoplus_{\gamma \in X_{#1}\\ \delta\ssupface\gamma} \bigwedge^{#1}\TT^\dual\gamma\otimes H^{#2}(\delta)}
\renewcommand{\bot}[2]{\displaystyle\bigoplus_{\gamma \in X_{#1}} \bigwedge^{#1}\TT^\dual\gamma\otimes H^{#2}(\gamma)}
\begin{tikzpicture}[scale=.9]
\node (00) at (0,0) {$\displaystyle\bigoplus_{v\in X_0}H^{2p}(v)$};
\node[left = .5cm of 00] (0-1) {$\displaystyle\bigoplus_{v\in X_\f \\ \dims{v}=0}H^{2p}(v)$};
\node[right = of 00] (02) {$\top{0}{2p-2}$};
\node[right = of 02] (04) {$\top{1}{2p-4}$};
\node[below = of 00] (10) {};
\node (11) at ($(10) + .5*(00)+.5*(02)$) {$\bot{0}{2p}$};
\node (13) at ($(10) + .5*(02)+.5*(04)$) {$\bot{1}{2p-2}$};
\node[inner sep=10pt, right = 1.7cm of 13] (15) {$\cdots$};
\foreach \i/\j in {0/1, 2/1, 2/3, 4/3, 4/5} {
  \draw[->, shorten <= -20pt] (0\i) -- (1\j);
}
\draw[{Hooks[right, length=4pt]}->] (0-1) -- (00);
\end{tikzpicture}\hspace{3cm}\\[1cm]
\hspace{2cm}\begin{tikzpicture}
\node[inner sep=15pt] (00) at (0,0) {};
\node[above] (00') {$\cdots$};
\node[right = of 00] (02) {$\top{p-1}{0}$};
\node[right = of 02] (04) {$\displaystyle\bigoplus_{\gamma \in X_p}\bigwedge^p \TT^\dual \gamma$};
\node[anchor=base, right = of 04] (06) {$\displaystyle C^{p,p}_\trop(X)$};
\node[below left = of 02] (10) {};
\node[left = of 10] (1-1) {$\cdots$};
\node (11) at ($(10) + .5*(00)+.5*(02)$) {$\bot{p-1}{2}$};
\node[right = of 11] (13) {$\bot{p}{0}$};
\foreach \i/\j in {2/1} {
  \draw[->, shorten <= -20pt] (0\i) -- (1\j);
}
\draw[->, shorten <= -10pt] (00) -- (11);
\draw[->, shorten <= -15pt] (02.-12) -- (13);
\draw[->, shorten <= -10pt] (04.south) -- (13);
\draw[<-] (04) -- (06);
\end{tikzpicture}
\caption{Zigzag diagram giving the isomorphism $P^{0,2p} \simeq PH^{p,p}_\trop(X)$ \label{fig:zigzag}}
\end{figure}
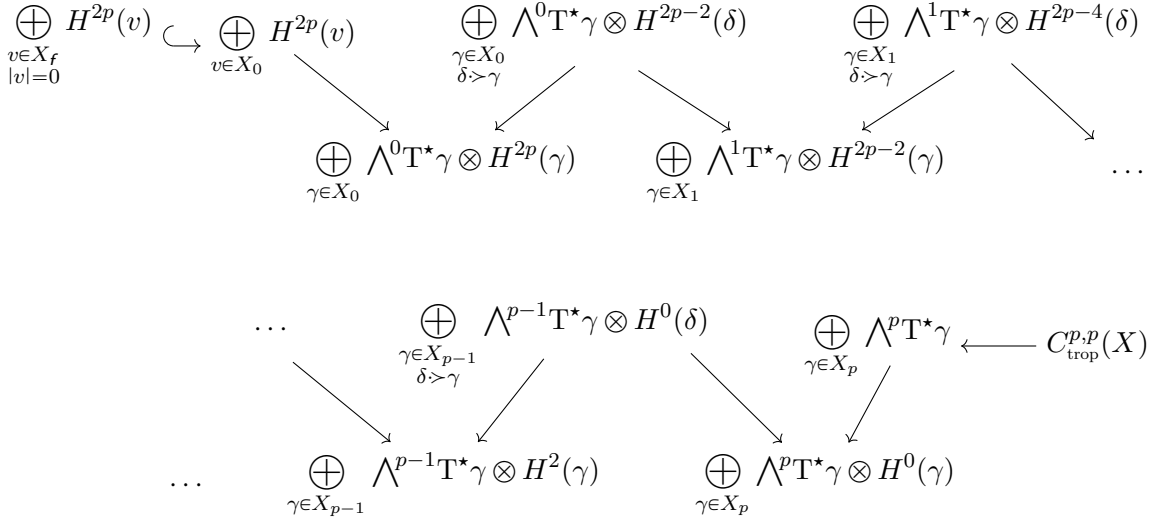

\smallskip
As already mentioned above, from the analysis carried out in Section 5 of \cite{AP-tht}, we infer that this isomorphism is given by the diagram of Figure \ref{fig:zigzag}.

In this diagram, the notation $\delta \ssupface \gamma$ in the direct sum running over faces $\gamma$ in $X$ of dimension $k$ means we consider any face $\delta$ of $X$ such that $\gamma$ is a face of codimension one in $\delta$ (and so $\delta$ has dimension $k+1$), and such that in addition, we have $\sed(\delta)=\sed(\gamma)$. If $\gamma \ssubface \delta$, we have two natural maps. The first one is
\begin{align*}
\cdot\wedge\nu^\dual_{\delta/\gamma}\colon \bigwedge^k\TT^\dual\gamma &\longrightarrow \bigwedge^{k+1}\TT^\dual\delta, \\
\alpha &\longmapsto \~\alpha\wedge\nu^\dual_{\delta/\gamma}.
\end{align*}
Here, $\~\alpha$ is any extension of $\alpha$ to $\bigwedge^k \TT^\dual\delta$, and we recall that $\nu^\dual_{\delta/\gamma}\in\TT^\dual_\delta$ is zero on $\TT\gamma$, nonnegative on $\TT\delta$ and induces a primitive vector in $M_\delta=N_\delta^\vee$. The image of $\alpha$ does not depend on the chosen extension. For the second map, recall first that by the Hodge isomorphism theorem, we have $H^{2k}(\delta) \simeq A^k(\Sigma^\delta)$. Moreover, any ray $\varrho$ in the star fan $\Sigma^\delta$ can be naturally identified with a ray $\varrho'$ of $\Sigma^\gamma$. Hence, we naturally get a map $\iota\colon A^k(\Sigma^\delta) \to A^k(\Sigma^\gamma)$ by mapping $\x_\varrho$ to $\x_{\varrho'}$. The Gysin map is then given by
\begin{align*}
\gys\colon A^k(\Sigma^\delta) &\longrightarrow A^{k+1}(\Sigma^\gamma), \\
x &\longmapsto \iota(x)\cdot \x_{\rho_{\delta/\gamma}},
\end{align*}
where $\rho_{\delta/\gamma}$ is the ray corresponding to $\delta$ is $\Sigma^\gamma$.

\medskip

The maps in Figure \ref{fig:zigzag} are naturally given by (the tensor product of) the two maps we just described, up to some signs. More precisely, the first two and the last two maps are the natural ones. For the arrows going to the left, if we fix a pair of faces $\delta \ssupface \gamma$ and a face $\gamma'$ of dimensions $k+1$, $k$, and $k$, respectively, then the map
\[ \bigwedge^k\TT^\dual\gamma \otimes H^{2(p-k-1)}(\delta) \to \bigwedge^k\TT^\dual\gamma' \otimes H^{2(p-k)}(\gamma') \]
is given by $0$ if $\gamma \neq \gamma'$, and by
\[ \sign(\gamma, \delta) \id\otimes\gys_{\delta\supface\gamma}, \]
if $\gamma=\gamma'$. For the arrows going to the right, and for a pair of faces $\delta \ssupface \gamma$ and a face $\delta'$ of dimensions $k+1$, $k$, and $k+1$, respectively, the map
\[ \bigwedge^k\TT^\dual\gamma \otimes H^{2(p-k-1)}(\delta) \to \bigwedge^k\TT^\dual\gamma' \otimes H^{2(p-k)}(\gamma') \]
is given by $0$ if $\delta\neq\delta'$, and by
\[ (\,\cdot\wedge\nu_{\delta/\gamma}^\dual) \otimes \id \]
if $\delta=\delta'$.

Note that all the maps except the first one are surjective. Therefore, to an element of $\bigoplus_{v \in X_{\f,0}} H^{2p}(v)$ we can associate an element of $C^{p,p}_\trop(X)$ by following the diagram, and making some choices. However, as explained above, the induced map on the level of the cohomology groups does not depend on these choices and is thus well-defined. The resulting map we obtain in this way is the restriction of the Steenbrink-Tropical comparison isomorphism of Theorem \ref{thm:steenbrink-tropical} which goes from $P^{0,2p}$ to $PH^{p,p}_\trop$.

\subsection{Pairing with Minkowski weights} \label{sec:pairing_with_mw} Let $\Y$ be the open part of $\X$ and let $Y$ be the unimodular triangulation induced by $X$ on $\Y$. Let $p\geq0$ and let $w$ be a Minkowski weight in $\MW_p(Y)$. The weight $w$ naturally induces a map $w\colon X_p \to \Q$ which we also denote by $w$.

\medskip

Let $v$ be a vertex on $X_\f$. Let $\alpha_v \in H^{2p}(v)$. We can represent $\alpha_v$ as follows
\[ \alpha_v = \sum_{\eta\supface v \\ \dims\eta = p} a_{v,\eta}\x_\eta \in A^p(\Sigma^v). \]
The Minkowski weight $w$ naturally induces an element $w_v$ in $\MW_p(\Sigma^v)$. Via the pairing $A^p(\Sigma^v) \times \MW_p(\Sigma^v) \to \Q$, we thus get the rational number
\[ \langle\alpha_v, w_v\rangle = \sum_{\eta\subface v \\ \dims\eta = p} a_{v,\eta}w_v(\eta). \]
The fact that this sum does not depend on the chosen representative of $\alpha_v$ is a consequence of the balancing condition (cf. \cite{AHK}). By summing over all the vertices of $X_\f$, we get a pairing between $\ker(N\colon \ST_1^{0,2p} \to \ST_1^{2,2p-2}) \simeq \bigoplus_v H^{2p}(v)$ and $\MW_p(Y)$. One can prove that this induces a pairing
\[ P^{0,2p} \times \MW_p(Y) \to \Q. \]

\smallskip
On the other hand, let $c \in C^{p,p}(X)$. We write $c = \sum_{\eta \in X_p} c_\eta$ with $c_\eta \in \SF^p(\eta)$. Then we get a pairing between $C^{p,p}(X)$ and $\MW_p(Y)$ by setting
\[ \langle c, w\rangle = \sum_{\eta \in X_p} w(\eta)c_\eta(\nvect_\eta), \]
where we recall that $\nvect_\eta$ is the canonical unit multivector in $\bigwedge^p\TT\eta$. Once again, this induces a pairing
\[ PH^{p,p}_\trop(X) \times \MW_p(Y) \to \Q. \]
Let $\alpha\in PH^{p,p}_\trop(X)$. Let $\comp C=(X_{(d-p)}, w)$ be the tropical cycle associated to $w$. Then the above pairing is compatible with the integration on tropical Dolbeault cohomology in the sense of \cite{JSS19}:
\[ \int_{\comp C}\alpha = \langle \alpha, w \rangle, \]
for any $\alpha \in PH^{p,p}(X) \hookrightarrow H^{p,p}_\Dolb(X, \R)$.

\smallskip
The following theorem states that the two pairing defined above are compatible with the Steenbrink-Tropical comparison theorem.

\begin{thm} \label{thm:mw_pairing_commute}
The restriction of the isomorphism of Theorem {\rm\ref{thm:steenbrink-tropical}} to $P^{0,2p} \simto PH^{p,p}_\trop(X)$ commutes with the pairing with $\MW_p(Y)$ in the sense that the following diagram commutes.
\[ \begin{tikzcd}[column sep=0]
P^{0,2p} \ar[rr, "\sim"] \ar[rd] && \ar[ld] PH^{p,p}_\trop(X) \\
& \MW_p(Y)^\dual &
\end{tikzcd} \]
\end{thm}

\begin{proof}
It suffices to extend the pairing in a consistent way to all the pieces appearing in the zigzag of Figure \ref{fig:zigzag}. We do it as follows. Let $w$ be a Minkowski weight in $\MW_p(Y)$, and let $w\colon X_p \to \Q$ be the corresponding map on $X_p$, as above. Let $\delta \ssupface \gamma$ be a pair of faces of dimensions $k+1$ and $k$. Let $\beta_\gamma \in \bigwedge^k\TT^\dual\gamma$ and let $\alpha_\delta\in H^{2(p-k-1)}(\delta)$. We set
\[ \langle \beta_\gamma \otimes \alpha_\delta, w \rangle := \sign(\gamma, \delta) \beta_\gamma(\nvect_\gamma) \langle\alpha_\delta, w_\delta\rangle \]
where $w_\delta$ is the natural element of $\MW_{d-k-1}(\Sigma^\delta)$ induced by $w$. In the same way, if $\alpha_\gamma \in H^{2(p-k)}(\gamma)$, then we set
\[ \langle \beta_\gamma \otimes \alpha_\gamma, w \rangle := \beta_\gamma(\nvect_\gamma) \langle\alpha_\gamma, w_\gamma\rangle. \]
A direct computation shows that these pairings commute with the maps of the zigzag. This concludes the proof of the theorem.
\end{proof}

\section{Tropical surviving and relative cohomology groups}\label{sec:surviving_relative}

In this section we define the tropical surviving and relative cohomology groups.

\medskip

First, we denote by $K^{\bul,\bul}$ the kernel of the monodromy operator $N\colon \ST_1^{\bul,\bul} \to \ST_1^{\bul,\bul}$. It is explicitly given by
\[ K^{a,b} := \ST_1^{a,b,a} = \bigoplus_{\delta\in X_\f \\ \dims\delta = a}H^b(\delta). \]
In particular, since the monodromy operator is injective on $\ST_1^{a,b}$ for $a<0$, the groups $K^{a,b}$ are all trivial for any value of $a<0$. We can naturally endow $K^{a,b}$ with the restriction of the differential $\d$, \ie, with $\i^*$. We set
\[ H_\s^{p,q}(X):=H^{q-p}(K^{\bul,2p}) \]
and we call it the \emph{surviving cohomology of the tropical variety $\X$ with respect to the triangulation $X$}.

\begin{remark} To justify the terminology, we note that in the approximable case, when the tropical variety appears as the tropical limit of a projective family of complex varieties, using the correspondence between faces of the tropical limit and the intersection of components of the special fiber, the surviving cohomology gives a way to compute the weight-graded pieces of the cohomology of $\fX$. More precisely, using the Deligne spectral sequence, the cohomology group $H_\s^{p,q}(\X)$ coincides with the weight $2p$ graded piece of the cohomology group $H^{p+q}$ of the special fiber. In the general case, a triangulation equipped with the double complex $K^{\bul, \bul}$ given above plays the role of a virtual special fiber.
\end{remark}

\medskip

Similarly, we denote by $R^{\bul,\bul}$ the cokernel of the monodromy operator $N\colon \ST_1^{\bul,\bul} \to \ST_1^{\bul,\bul}$. It is explicitly given by
\[ R^{a,b} := \ST_1^{a,b,-a} = \bigoplus_{\delta\in X_\f \\ \dims\delta = -a}H^{2a+b}(\delta). \]
Since the monodromy operator is surjective on $\ST_1^{a,b}$ for any $a>0$, the group $R^{a,b}$ is trivial for $a>0$. We can endow $R^{a,b}$ with the restriction of the differential $\d$, \ie, with $\gys$. We set
\[ H_\rel^{p,q}(X):=H^{q-p}(R^{\bul,2p}) \]
and we call it the \emph{relative cohomology of $X$ relative to the triangulation $X$}.

\medskip

To justify the name of the cohomology $H_\rel^{\bul,\bul}(\X)$, one can compare the proposition which follows with the two triangles in the setting of complex manifolds, that we reproduce here from Section \ref{sec:classical}.
\[ \begin{tikzcd}[column sep=tiny]
H^\bul(\ufX) \ar[rr, "N"] && H^\bul(\ufX) \ar[dl, "+1"] \\
 & \ar[ul] H^\bul(\fX^*) &
\end{tikzcd} \qquad \begin{tikzcd}[column sep=tiny]
H^\bul(\fX) \ar[rr] && H^\bul(\fX^*) \ar[dl, "+1"] \\
 & \ar[ul] H^\bul(\fX, \fX^*) &
\end{tikzcd}
\]

If $\Phi\colon A^\bul \to B^\bul$ is a morphism of complexes, we denote by $\Cone^\bul(\Phi)$ the mapping cone of $\Phi$. We recall that it is defined by the total complex of the following double complex
\[ \begin{tikzcd}
  \cdots \rar& A^{0} \rar\dar["\Phi"]& \underset{\bullet}{A^{1}} \rar\dar["\Phi"]& A^{2} \rar\dar["\Phi"]& \cdots \\
  \cdots \rar& \underset{\bullet}{B^{0}} \rar& B^{1} \rar& B^{2} \rar& \cdots
\end{tikzcd} \]
where the pieces which form the degree zero part of the total complex are indicated by a dot. Note that the differentials on the first row are minus the differentials on $A^\bul$.

\begin{prop}
For any even integer $p$, we have the two following distinguished triangles in the derived category of bounded cochain complexes
\[ \begin{tikzcd}[column sep=tiny]
\ST_1^{\bul,p+2} \ar[rr, "N"] && \ST_1^{\bul,p}[2] \ar[dl, "+1"] \\
 & \ar[ul] \Cone^\bul(N)[-1]
\end{tikzcd} \qquad \begin{tikzcd}[column sep=tiny]
K^{\bul,p+2} \ar[rr] && \Cone^\bul(N)[1] \ar[dl, "+1"] \\
 & \ar[ul, "0"'] R^{\bul,p} &
\end{tikzcd} \]
in which $\Cone^\bul(N)$ refers more precisely to \[ \Cone^\bul\Bigl(N\colon \ST_1^{\bul,p+2} \to \ST_1^{\bul,p}[2]\Bigr). \]
\end{prop}

\begin{proof}
The first triangle is a distinguished triangle by definition.

The shifted mapping cone
\[ \Cone\Bigl(N\colon \ST_1^{\bul,p+2} \to \ST_1^{\bul,p}[2]\Bigr)[-1] \]
is given by the total complex associated to the following double complex
\[ \begin{tikzcd}
  \cdots \rar["\d"] & \ST_1^{-3,p+2} \rar["\d"] \dar["N"] & \ST_1^{-2,p+2} \rar["\d"] \dar["N"] & \ST_1^{-1,p+2} \rar["\d"] \dar["N"] & \underset{\bullet}{\ST_1^{0,p+2}} \rar["\d"] \dar["N"] & \cdots \\
  \cdots \rar["\d"] & \ST_1^{-1,p} \rar["\d"] & \ST_1^{0,p} \rar["\d"] & \underset{\bullet}{\ST_1^{1,p}} \rar["\d"] & \ST_1^{2,p} \rar["\d"] & \cdots
\end{tikzcd} \]
Recall that $R^{\bul,p}$ is the cokernel of $N\colon \ST_1^{\bul,p+2}[-2]\to \ST_1^{\bul,p}$. Thus, we can naturally project the second row onto $R^{\bul,p}[1]$ to get the map on the right hand side of the second triangle.

\smallskip
The shifted mapping cone
\[ \Cone^\bul\Bigl(K^{\bul,p+2} \to \Cone^\bul(N)[1]\Bigr)[-1] \]
is given by the total complex associated to the following double complex
\[ \begin{tikzcd}
  \cdots \rar &  0                               \rar \dar[hook, "\iota"] &  0                                  \rar \dar[hook, "\iota"] &  \underset{\bullet}{K^{0,p+2}}  \rar \dar[hook, "\iota"] &  K^{1,p+2}     \rar \dar[hook, "\iota"] &  \cdots \\
  \cdots \rar &  \ST_1^{-2,p+2}                   \rar \dar["N"] &  \underset{\bullet}{\ST_1^{-1,p+2}}  \rar \dar["N"] &  \ST_1^{0,p+2}                   \rar \dar["N"] &  \ST_1^{1,p+2}  \rar \dar["N"] &  \cdots \\
  \cdots \rar &  \underset{\bullet}{\ST_1^{0,p}}            \rar &  \ST_1^{1,p}                                   \rar &  \ST_1^{2,p}                               \rar &  \ST_1^{3,p}              \rar &  \cdots
\end{tikzcd} \]
where the dots indicate the pieces of degree zero in the total complex. To simplify the notation, we denote this total complex by $T^\bul$.

\medskip

We just have to prove that the natural projection $\pi$ of the third row to $R^{\bul,p}$ induces a quasi-isomorphism between $T^\bul$ and $R^{\bul,p}$.

First, we can compute the cohomology of $T^\bul$ by the spectral sequence given by the filtration by rows. For the first page, we simply obtain $R^{\bul,p}$ (on the third row but without shift of the degrees). Therefore, $H^\bul(T^\bul) \simeq H^\bul(R^{\bul,p})$. It remains to prove that this isomorphism is given by the projection $\pi$.

\smallskip
To do so, it suffices to prove that for any cocycle $x$ of $T^a$, $a\leq0$, if $\pi(x)$ is a coboundary, then $x$ is a coboundary.

Let $x$ be such an element. In what follows, we denote by $\d$ the horizontal coboundary maps and by $\d_T$ the coboundary maps of the total complex $T^\bul$.

Write $x=x_1+x_2+x_3\in T^a$ with $x_1\in\ST_1^{a,p}$, $x_2\in\ST_1^{a-1,p+2}$ and $x_3\in K^{a,p+2}$. Let $y'\in R^{a-1,p}$ such that $\d y'=\pi(x)$. Let $y_1\in\ST_1^{a-1,p}$ such that $\pi(y_1)=y'$. Then $\pi(x_1-\d y_1)=0$ and there exists $y_2\in \ST_1^{a-2,p+2}$ such that $Ny_2 = x_1-\d y_1$. Set $y=y_1+y_2$. Since $x$ is a cocycle, $\d_T (x-\d_T y)=0$. In particular, $N(x_2-\d y_2)=-\d(x_1-\d y_1-Ny_2)=0$. Since $N\colon \ST^{a-1,p+2}\to\ST^{a+1,p}$ is injective for $a\leq0$, we deduce that $x_2-\d y_2=0$. Thus, $\iota x_3=-\d(x_2-\d y_2)=0$. By the injectivity of $\iota$, $x_3=0$. Finally, $x=\d_T(y_1+y_2)$ is a coboundary, which concludes the proof.
\end{proof}

\section{Tropical Clemens-Schmid sequence}
\label{sec:CS_sequence}

The aim of this section is to prove Theorem~\ref{thm:CS}.

\subsection{Main theorem in an abstract setting}\label{sec:cs_abstract} Let $(C^\bul, \d_C)$ and $(D^\bul, \d_D)$ be two bounded cochain complexes. Let $L\colon C^\bul\to D^\bul[2]$ be a morphism of cochain complexes. We say that $L$ verifies the Hard Lefschetz property around $0$ that we denote again by $\HL$, if for any nonnegative integer $k$, the map $C^{k-1} \xrightarrow{L} D^{k+1}$ is injective provided $k\leq0$ and surjective provided $k\geq0$.

\smallskip

Let $H^\bul(C^\bul)$ and $H^\bul(D^\bul)$ be the cohomology of $C^\bul$ and $D^\bul$, respectively. Assume that $\d_DL=L\d_C$, and that both $L\colon C^\bul\to D^\bul[2]$ and the induced map on cohomology $L\colon H^\bul(C^\bul)\to H^\bul(D^\bul)[2]$ verify the hard Lefschetz property around $0$ $\HL$. Finally, set
\[ K^\bul:=\ker\bigl(L\colon C^\bul\to D^\bul[2]\bigr) \qquad \textrm{and} \qquad R^\bul:=\coker\bigl(L\colon C^\bul[-2]\to D^\bul\bigr). \]

Then we can endow $K^\bul$ and $R^\bul$ with the corresponding differential operators to get two cochain complexes. We denote their respective cohomology groups by $H^\bul(K^\bul)$ and $H^\bul(R^\bul)$.

\begin{thm}[Abstract Clemens-Schmid sequence] \label{thm:main}
With the above notations, we get the two following long exact sequences depending on the parity of $k$
\[ \cdots\to H^{k-1}(K^\bul) \to H^{k-1}(C^\bul) \xrightarrow{L} H^{k+1}(D^\bul) \to H^{k+1}(R^\bul) \to H^{k+1}(K^\bul) \to \cdots \]
\end{thm}

\begin{proof}
First, by the definitions which preceded, we have an exact sequence
\[ 0 \to K^\bul[-1] \to C^\bul[-1] \xrightarrow{L} D^\bul[1] \to R^\bul[1] \to 0. \]
By the injectivity of $L$ on negative degrees, $K^k=0$ if $k<0$. In the same way, by the surjectivity of $L$ on degrees at least $-1$, we get $R^k=0$ if $k>0$. We can thus unfold the previous exact sequence of complexes into the diagram of Figure \ref{fig:diagram_proof_CS} where rows are all exact.
\begin{figure}[ht]
\newcommand{\ear}{\ar[r, very thin, -latex]}
\[ \begin{tikzcd}[column sep=scriptsize]
          &                   & \cdvdots     \dar & \cdvdots     \dar & \cdvdots     \dar &          \\
          & 0        \ear     & C^{-3}   \ear\dar & D^{-1}   \ear\dar & R^{-1}   \ear\dar & 0        \\
          & 0        \ear     & C^{-2}   \ear\dar & D^0      \ear\dar & R^0      \ear\dar & 0        \\
          & 0        \ear\dar & C^{-1}   \ear\dar & D^1      \ear\dar & 0                 &          \\
  0  \ear & K^0      \ear\dar & C^0      \ear\dar & D^2      \ear\dar & 0                 &          \\
  0  \ear & K^1      \ear\dar & C^1      \ear\dar & D^3      \ear\dar & 0                 &          \\
          & \cdvdots          & \cdvdots          & \cdvdots          &                   &          \\
\end{tikzcd} \]
\caption{\label{fig:diagram_proof_CS}}
\end{figure}

For a cochain complex $E^\bul$ and an integer $l$, denote by $E_{\leq l}^\bul$ the truncation of $E^\bul$ to parts of degree at most $l$. From the above discussion, we thus get a short exact sequence \[ 0 \to C^\bul_{\leq-1}[-1] \xrightarrow{L} D^\bul_{\leq1}[1] \to R^\bul[1] \to 0, \]
given by the rows on top of the diagram in Figure \ref{fig:diagram_proof_CS} up to the middle row.

This short exact sequence induces a long exact sequence
\begin{equation} \label{eqn:long_seq_1}
\begin{aligned}
\cdots \to H^k(C^\bul) \xrightarrow{L} H^{k+2}(D^\bul) \to H^{k+2}(R^\bul) \to H^{k+1}(C^\bul) \to \cdots \qquad\qquad & \\
\cdots \to H^{-2}(C^\bul) \xrightarrow{L} H^0(D^\bul) \to H^0(R^\bul) \to\ ?&
\end{aligned}
\end{equation}
where the question mark on the right hand side is the cokernel of the map $C^{-2}\xrightarrow{\d}C^{-1}$, which is a priori different from $H^{-1}(C^\bul)$.

\medskip

Actually we can split this long exact sequence into short exact sequences. Indeed, we know that $H^{k-1}(C^\bul)\xrightarrow{L} H^{k+1}(D^\bul)$ is injective if $k\leq0$. Thus, maps of the form $H^k(R^\bul) \to H^{k-1}(C^\bul)$ are zero in the sequence for any $k\leq-1$. For such an integer $k$, we get a short exact sequence
\[ 0 \to H^{k-2}(C^\bul) \xrightarrow{L} H^k(D^\bul) \to H^k(R^\bul) \to 0. \]
By a symmetric argument, we get short exact sequences
\[ 0 \to H^k(K^\bul) \to H^k(C^\bul) \xrightarrow{L} H^{k+2}(D^\bul) \to 0, \]
for any $k\geq1$. By \HL, we also know that we have an isomorphism
\[ 0 \to H^{-1}(C^\bul) \to H^1(D^\bul) \to 0. \]

\medskip

Gluing all these short exact sequences, we almost get the long exact sequences of the theorem. In fact, we directly get the long exact sequence in which all the degrees are odd integers, \ie, with $k$ in the statement of the theorem is even. To see this, note that for a positive even integer $k$, we have
\begin{gather*}
\cdots \to \underbrace{H^{-k-1}(K^\bul)}_0 \to H^{-k-1}(C^\bul) \xrightarrow{L} H^{-k+1}(D^\bul) \to H^{-k+1}(R^\bul) \to \underbrace{H^{-k+1}(K^\bul)}_0 \to \cdots \\
\cdots \to \underbrace{H^{-1}(K^\bul)}_0 \to H^{-1}(C^\bul) \xrightarrow{L} H^{1}(D^\bul) \to \underbrace{H^1(R^\bul)}_0 \to \cdots \\
\cdots \to \underbrace{H^{k-1}(R^\bul)}_0 \to H^{k-1}(K^\bul) \to H^{k-1}(C^\bul) \xrightarrow{L} H^{k+1}(D^\bul) \to \underbrace{H^{k+1}(R^\bul)}_0 \to H^{k-1}(K^\bul) \to \cdots
\end{gather*}
which is exactly the above exact sequences, combined together.

For the other exact sequence in the theorem, \ie, when all the degrees are even, we can apply a similar argument as above to treats all the other cases and reduce to proving the exactness of the following six-term sequence
\begin{equation} \label{eq:six-term-exact-sequence}
0 \to H^{-2}(C^\bul) \to H^0(D^\bul) \xrightarrow{\d^{-1}} H^0(R^\bul) \xrightarrow{\d^0} H^0(K^\bul) \xrightarrow{\d^1} H^0(C^\bul) \to H^2(D^\bul) \to 0.
\end{equation}

The exactness of the beginning of this sequence is a consequence of \eqref{eqn:long_seq_1} and the injectivity of $L\colon H^{-2}(C^\bul) \to H^0(D^\bul)$. By a symmetric argument, we infer the exactness of the end of the sequence. It thus remains to describe the central map $\d^0$, and to prove the exactness of the sequence at other places, \ie, to show that $\Im(\d^{-1})=\ker(\d^0)$ and $\Im(\d^0)=\ker(\d^1)$.

\medskip

The end of the proof is essentially a diagram chasing. The definition of $\d^0$ is given by the diagram depicted in Figure~\ref{fig:definition-d0}.

In this diagram, the usual maps have been removed for the sake of an increase in readability. The drawn arrows indicate in which order the diagram chasing is done. For instance, the red part has to be read as follows: take an element $r$ in $R^0$. There exists a preimage $c$ of $r$ in $D^0$, etc. An exclamation mark above a red arrow indicates that the preimage is unique.

The red part defines $\d^0$ by mapping $r$ onto $a''$. Moreover, it shows that $a''$ belongs to $\ker(K^0\to K^1)$. The blue part proves that $a''$ is independent of the choice of $c$: if instead of $c$ one chooses $\~c$, then $b''$ will not change. Finally, the green part proves that if $r$ is a coboundary, then $c'=0$, which implies that $a''=0$. Altogether, we thus get a well-defined map $\d^0\colon H^0(R^\bul)\to H^0(K^\bul)$.

\begin{figure}
\begin{tikzpicture}[commutative diagrams/every diagram]
\matrix[matrix of math nodes, name=m, commutative diagrams/every cell, gray!50] {
     & 0    & C^{-3}  & D^{-1} & R^{-1} & 0  \\
     & 0    & C^{-2}  & D^0    & R^0    & 0  \\
     & 0    & C^{-1}  & D^1    & 0      &    \\
  0  & K^0  & C^0     & D^2    & 0 \\
  0  & K^1  & C^1     & D^3    & 0 \\
};
\foreach \i/\j in {1/2, 1/6, 2/2, 2/6, 3/2, 3/5, 4/1, 4/5, 5/1, 5/5} {
  \path (m-\i-\j) node {$0$};
}
\foreach \i/\j/\k in {1/4/-c, 1/5/-r, 2/3/b, 2/4/c, 2/5/r, 3/3/b', 3/4/c', 4/2/a'', 4/3/b'', 4/4/c'', 5/3/b''', 5/2/a'''} {
  \node[below right, white] (\k) at (m-\i-\j) {$\k$};
  \node[below right, white] (2\k) at (\k) {$\k$};
}
\draw[red!70!black, commutative diagrams/.cd, every arrow]
  (r) node {$r$}
    edge["$\exists$"'] (c) (c) node {$c$} (c)
    edge (c') (c') node {$c'$}
    edge (c'') (c'') node {$0$}
  (c')
    edge["$!$"'] (b') (b') node {$b'$}
    edge (b'') (b'') node {$b''$}
    edge["$!$"'] (a'') (a'') node {$a''$}
  (b'')
    edge (b''') (b''') node {$0$}
    edge["$!$"'] (a''') (a''') node {$0$}
    edge (a'');
\scriptsize
\draw[green!70!black, commutative diagrams/.cd, every arrow]
  (r)
    edge (-r) (-r) node {$\varrho$}
    edge["$\exists$"] (-c) (-c) node {$\gamma$}
    edge (c)
  (c') node[right] (3c') {${}=0$}
  (c)
    edge (3c');
\draw[blue!70!black, commutative diagrams/.cd, every arrow]
  (2c) node[fill=white, fill opacity=.6, text opacity=1, inner sep=0pt, outer sep=2pt] (2c) {$c\!-\!\~c$}
    edge (2r) (2r) node {$0$}
  (2c)
    edge (b) (b) node {$\beta$}
    edge (2b') (2b') node {$\beta'$}
    edge (2b'') (2b'') node {$0$};
\end{tikzpicture}
\caption{Definition of $\d^0$ in the six-term exact sequence~\eqref{eq:six-term-exact-sequence}.}
\label{fig:definition-d0}
\end{figure}

Moreover, one can deduce from this diagram that the decomposition $\d^0\d^{-1}$ is zero. Indeed, if $c\in\ker(D^0\to D^1)$ and if we denote by $\class(c)$ its class in $H^0$, then we have $\d^0\d^{-1}\bigl(\class(c)\bigr)=\class{(a'')}$. But $c'=0$, and so $a''=0$.

By a somehow similar reasoning, we get $\d^1\d^0=0$. Indeed, $\d^1\d^0\bigl(\class{(r)}\bigr)=\class{(b'')}$ and $b''$ is a coboundary.

\medskip

It remains to prove that $\ker(\d^0)\subseteq\Im(\d^{-1})$ and that $\ker(\d^1)\subseteq\Im(\d^0)$. For the first inclusion, we use the diagram depicted in Figure~\ref{fig:first-injectivity}. The red elements are defined as previously. If $\d^0\bigl(\class{(r)}\big)=0$, then $a''=0$, and the blue part of the diagram follows. Then $b'$ is in the kernel of the map $C^{-1}\to C^0$, and the map $L\colon H^{-1}(C^\bul)\simto H^1(D^\bul)$ is an isomorphism mapping $\class{(b')}$ onto $\class{(c')}$. Since $c'$ is a coboundary, $\class{(c')}=0$. Hence $\class{(b')}=0$ and $b'$ is a coboundary. From this, we get the green part. Finally, $c-\~c\in\ker(D^0\to D^1)$ and $\d^{-1}\class{(c-\~c)}=\class(r)$.

\begin{figure}[!ht]
\begin{tikzpicture}[commutative diagrams/every diagram]
\matrix[matrix of math nodes, name=m, commutative diagrams/every cell, gray!50] {
     & 0    & C^{-2}  & D^0    & R^0    & 0  \\
     & 0    & C^{-1}  & D^1    & 0      &    \\
  0  & K^0  & C^0     & D^2    & 0 \\
};
\foreach \i/\j in {1/2, 1/6, 2/2, 2/5, 3/1, 3/5} {
  \path (m-\i-\j) node {$0$};
}
\foreach \i/\j/\k in {1/3/b, 1/4/c, 1/5/r, 2/3/b', 2/4/c', 3/2/a'', 3/3/b'', 3/4/c''} {
  \node[below right, white] (\k) at (m-\i-\j) {$\k$};
}
\draw[red!70!black]
  (a'') node {$a''$}
  (b'') node {$b''$}
  (c'') node {$0$}
  (b') node {$b'$}
  (c') node {$c'$}
  (c) node {$c$}
  (r) node {$r$};
\draw[blue!70!black, commutative diagrams/.cd, every arrow]
  (a'') node[right] (2a'') {${}=0$}
  (b'') node[right] (2b'') {${}=0$}
  (2a'')
    edge (b'')
  (b')
    edge["$\sim$"] (c');
\draw[green!70!black, commutative diagrams/.cd, every arrow]
  (c) node[below right] (2c) {$\~c$}
  (r) node[below right] (2r) {$0$}
  (b')
    edge["$\exists$"] (b) (b) node {$b$}
  (b)
    edge (2c) (2c)
    edge (c')
  (2c)
    edge (2r);
\end{tikzpicture}
\caption{Proof of the inclusion $\ker(\d^0)\subseteq\Im(\d^{-1})$.}
\label{fig:first-injectivity}
\end{figure}

For the second inclusion $\ker(\d^1)\subseteq\Im(\d^0)$, we use the diagram in Figure~\ref{fig:second-injectivity}. Here $a''$ is an element of the kernel $\ker(K^0\to K^1)$ such that $\d^1\class{(a'')}=0$. Thus, $b''$ is a coboundary, and the red part follows. For the blue part, notice that $c'$ is in the kernel of $D^1\to D^2$. We get a cocycle $\~{b'}\in C^{-1}$ such that $L\class{(\~{b'})}=\class{(c')}$. In particular, there exists a coboundary $\~{c'}\in D^1$ such that $L\~{b'}=c'-\~{c'}$. For the green part, we consider the element $b'-\~{b'}$, and every green arrow is clear. Finally, we construct a cocycle $r$ in $R^0$ such that $\d^0\class{(r)}=\class{(a'')}$.

\begin{figure}[!ht]
\begin{tikzpicture}[commutative diagrams/every diagram]
\matrix[matrix of math nodes, name=m, commutative diagrams/every cell, gray!50] {
     & 0    & C^{-2}  & D^0    & R^0    & 0  \\
     & 0    & C^{-1}  & D^1    & 0      &    \\
  0  & K^0  & C^0     & D^2    & 0 \\
};
\foreach \i/\j in {1/2, 1/6, 2/2, 2/5, 3/1, 3/5} {
  \path (m-\i-\j) node {$0$};
}
\foreach \i/\j/\k in {1/3/b, 1/4/c, 1/5/r, 2/3/b', 2/4/c', 3/2/a'', 3/3/b'', 3/4/c''} {
  \node[below right, white] (\k) at (m-\i-\j) {$\k$};
}
\draw[red!70!black, commutative diagrams/.cd, every arrow]
  ([shift={(0,-1)}]a'') node (a''') {$0$}
  (a'') node {$a''$}
    edge (a''')
  (a'')
    edge (b'') (b'') node {$b''$}
    edge["\tiny $\exists$"] (b') (b') node {$b'$}
    edge (c') (c') node {$c'$}
  (b'')
    to (c'') (c'') node {$0$};
\draw[blue!70!black, commutative diagrams/.cd, every arrow]
  (b') node[below right] (2b') {$\~{b'}$}
  (b'') node[below right] (2b'') {$0$}
  (c') node[below right] (2c') {$c'-\~{c'}$}
  (c')
    edge (2b') (2b')
    edge (2b'')
  (2b')
    to (2c');
\draw[green!70!black, commutative diagrams/.cd, every arrow]
  (m-2-3) node[above right, inner sep=1pt] (3b') {$b'-\~{b'}$}
  (m-2-4.30) node (3c') {$\~{c'}$}
  (3b'.south west)
    edge[semithick, bend right] (b'')
  (3b')
    edge (3c') (3c')
    edge["\tiny $\exists$", near end] (c) (c) node {$c$}
    to (r) (r) node {$r$};
\end{tikzpicture}
\caption{Proof of the inclusion $\ker(\d^1)\subseteq\Im(\d^0)$.}
\label{fig:second-injectivity}
\end{figure}

Combining all the previous results together, we get that $\d^0$ is well-defined, that $\Im(\d^{-1})=\ker(\d^0)$ and that $\Im(\d^0)=\ker(\d^1)$. This implies that the sequence
\[ 0 \to H^{-2}(C^\bul) \to H^0(D^\bul) \xrightarrow{\d^{-1}} H^0(R^\bul) \xrightarrow{\d^0} H^0(K^\bul) \xrightarrow{\d^1} H^0(C^\bul) \to H^2(D^\bul) \to 0 \]
is exact. Gluing all these exact sequences together, we get the long exact sequence of the theorem consisting only of even degrees. This means for an odd integer $k \geq 3$, we have the exact sequence

\begin{gather*}
\cdots \to \underbrace{H^{-k-1}(K^\bul)}_0 \to H^{-k-1}(C^\bul) \xrightarrow{L} H^{-k+1}(D^\bul) \to H^{-k+1}(R^\bul) \to \underbrace{H^{-k+1}(K^\bul)}_0 \to \cdots \\
\begin{aligned}
\cdots \to \underbrace{H^{-2}(K^\bul)}_0 \to H^{-2}(C^\bul) \to H^{0}(D^\bul) \to H^0(R^\bul) \to H^0(K^\bul) \to H^0(C^\bul)\hspace{2cm}& \\*[-1em]
\hspace{6cm} \to H^2(D^\bul) \to \underbrace{H^{2}(R^\bul)}_0 \to \cdots&
\end{aligned} \\
\cdots \to \underbrace{H^{k-1}(R^\bul)}_0 \to H^{k-1}(K^\bul) \to H^{k-1}(C^\bul) \xrightarrow{L} H^{k+1}(D^\bul) \to \underbrace{H^{k+1}(R^\bul)}_0 \to H^{k+1}(K^\bul) \to \cdots
\end{gather*}
This concludes the proof of the theorem.
\end{proof}

\subsection{Tropical Clemens-Schmid sequence} We now derive the tropical Clemens-Schmid sequence from the theorem established in the previous section. So let $\X$ be a rationally triangulable smooth projective tropical variety and fix a unimodular triangulation $X$ of $\X$. Consider the monodromy operator $N$, of bidegree $(2,-2)$, acting on the tropical Steenbrink double sequence $\ST_1^{\bul,\bul}$ associated to $N$. We know that $N$ verifies the Hard Lefschetz property $\HL$. In particular, $N\colon \ST_1^{a-1,b} \to \ST_1^{a+1,b-2}$ is injective if $a\leq0$ and surjective if $a\geq0$. Thus, the restriction $N\colon \ST_1^{\bul,b} \to \ST_1^{\bul,b-2}[2]$ is a monodromy operator verifying $\HL$ around $0$ in the sense of Section \ref{sec:cs_abstract}. Moreover, we know that $N$ induces an operator on $\ST_2^{\bul,\bul}=\ST_\infty^{\bul,\bul}$ which verifies $\HL$. Thus, $N\colon H^\bul\bigl(\ST_1^{\bul,b}\bigr) \to H^\bul\bigl(\ST_1^{\bul,b-2}\bigr)[2]$ also verifies $\HL$ around $0$.

\medskip

We can apply Theorem \ref{thm:main}, setting $p=b/2$ and $q=a+b/2$, to get the long exact sequence
\[ \cdots \to H_\s^{p,q}(X,\Q) \to H_\trop^{p,q}(\X,\Q) \xrightarrow{N} H_\trop^{p-1,q+1}(\X,\Q) \to H_\rel^{p-1,q+1}(X,\Q) \to H_\s^{p,q+2}(X,\Q) \to \cdots \]
Setting
\[ H^k(\X,\Q):=\bigoplus_{p+q=k} H_\trop^{p,q}(\X,\Q), \qquad \textrm{and}\]
\[\qquad H^k_\s(X,\Q) := \bigoplus_{p+q=k} H_\s^{p,q}(X,\Q), \qquad H^k_\rel(X, \Q) := \bigoplus_{p+q=k} H_\rel^{p,q}(X,\Q),\]
we can sum up the above exact sequences to get the tropical Clemens-Schmid exact sequence
\[ \cdots \to H_\s^k(X, \Q) \to H^k(\X,\Q) \xrightarrow{N} H^k(\X,\Q) \to H_\rel^k(X,\Q) \to H_\s^{k+2}(X,\Q) \to \cdots \]
as required. This finishes the proof of Theorem~\ref{thm:CS}. \qed

\section{Existence of cycles with a given Hodge class: proof of Theorem~\ref{thm:HC}}\label{sec:proof_HC}

In this section, we prove Theorem~\ref{thm:HC}. Let $\X$ be a smooth projective tropical variety of dimension $d$. Assume that $\X$ is rationally triangulable. We need to prove that any Hodge class $\alpha$, \ie, any element $\alpha$ in $\ker\bigl(N\colon H_\trop^{p,p}(\X, \Q) \to H_\trop^{p-1,p+1}(\X, \Q)\bigr)$, is represented by a tropical cycle. The statement that the image of the tropical cycle class map is in the kernel of the tropical monodromy follows from Theorem~\ref{thm:eigenwave_monodromy} and the same result proved by Mikhalkin and Zharkov~\cite{MZ14} for the eigenwave operator.

\smallskip
Replacing the underlying lattice by a rational multiple, we fix a unimodular triangulation $X$ of $\X$ with open part $Y$, and we will show the existence of a Minkowski weight in $\MW_{d-p}(Y)$ whose associated cycle in $X$ represents the Hodge class $\alpha$.

\medskip

By the tropical Clemens-Schmid exact sequence, we know that
\[ \ker\Bigl(N\colon H_\trop^{p,p}(\X, \Q) \to H_\trop^{p-1,p+1}(\X, \Q)\Bigr) = \Im\Bigl(H_\s^{p,p}(X,\Q)\to H_\trop^{p,p}(\X,\Q)\Bigr). \]
Moreover, by definition, we have
\[H_s^{p,p}(X,\Q) = \ker\Bigl( \bigoplus_{v \in X_{\f,0}} H^{2p}(v) \longrightarrow \bigoplus_{e \in X_{\f,1}} H^{2p}(e)\Bigr).\]
Thus, every Hodge class is represented by a cocycle of
\[ K^{0,2p}=\bigoplus_{v\in X_{\f,0}} H^{2p}(v). \]

The Hodge class $\alpha$ corresponds therefore to a collection of elements $\alpha_v\in H^{2p}(v)$, for $v$ a vertex in the finite part $X_\f$, with the compatibility condition
\[\forall\:e =uv \in X_{\f,1}, \qquad \i_{v\subface e}^*(\alpha_v) = \i_{u\subface e}^*(\alpha_u). \]

For each simplex $\delta$ in $Y$ and integer $k$, let $\MW_{k}(\delta):= \MW_{k}(\Sigma^\delta)$. By Theorems~\ref{thm:HI} and~\ref{thm:hc-local} we have $H^{2p}(v, \Q) \simeq A^p(\Sigma^v, \Q) \simeq \MW_{d-p}(v, \Q)$. The fact that $\alpha$ is a cocycle and the isomorphism $H^{2p}(e,\Q) \simeq \MW_{d-p}(e,\Q)$ now imply that the Minkowski weights around each vertex coincides on the star fans of the incident edges. This shows that $\alpha_v$ glue together and produce a Minkowski weight $C= (Y_{(d-p)}, w)$ of dimension $d-p$ on $Y$. As we explained in Section~\ref{sec:minkowski_weight}, this provides a tropical cycle $\comp C$ on $\X$.

It remains to prove that the cohomology class $\class(\comp C)$ in $H^{p,p}_\trop(\X, \Q)$ associated to $\comp C$ coincides with $\alpha$.
Let $\beta = (\beta_v)_v \in K^{0,2d-2p}$. For each vertex $v$ in $X_{\f,0}$ and for each face $\eta$ of dimension $d-p$ containing $v$, we choose a rational coefficient $b_{v,\eta} \in \Q$ such that
\[ \beta_v = \sum_{\eta \supface v} b_{v,\eta}\x_\eta \in A^{d-p}(\Sigma^v, \Q), \]
where the sum is over faces $\eta$ of dimension $d-p$. The local isomorphism between the Chow ring and the Minkowski weights given in Theorem~\ref{thm:hc-local} verifies the following equality for any $v\in X_{\f,0}$
\[ \deg(\alpha_v \cdot \beta_v) = \sum_{\eta \supface v} w(\eta)b_{v,\eta}. \]
Hence,
\[ \deg(\alpha\cdot\beta) = \sum_{v \in X_{\f,0}} \sum_{\eta \supface v} w(\eta)b_{v,\eta}. \]
Moreover, by Section \ref{sec:pairing_with_mw}, the pairing between $w$ and $\beta$ is
\[ \langle \beta, w\rangle = \sum_{v\in X_{\f,0}} \sum_{\eta \supface v} w(\eta)b_{v,\eta}. \]
Using Theorem \ref{thm:mw_pairing_commute}, we get
\[ \int_{\comp C}\beta = \langle \beta, w \rangle = \deg(\alpha \cdot \beta) \]
which shows that the class of $\comp C$ coincides with $\alpha$.

We proved that any Hodge class comes from a tropical cycle, which is the statement of Theorem~\ref{thm:HC}. \qed

\section{Proof of Theorem~\ref{thm:standard}}\label{sec:standard}

In this final section, we prove the equivalence of numerical and homological equivalence for tropical cycles on smooth projective tropical varieties which admit a rational triangulation.

Let $\X$ be such a tropical variety of dimension $d$ and consider a unimodular triangulation $X$ of $\X$, which exists after replacing the lattice by a rational multiple. The monodromy operator $N$ is a Lefschetz operator on $\ST_1^{\bul, \bul}$ and induces a Lefschetz operator on the cohomology of the tropical Steenbrink double sequence. From this, and using the Steenbrink-Tropical comparison theorem, we infer that $H^{p,p}_\trop(\X, \Q) \simeq H^0(\ST^{\bullet, 2p}(X))$ can be decomposed as a direct sum of the form
\[H^{p,p}(\X, \Q) = \ker\Bigl(N\colon H^{p,p}(\X) \to H^{p-1,p+1}(\X)\Bigr) \bigoplus \Im\Bigl(N\colon H^{p+1,p-1}(\X) \to H^{p,p}(\X)\Bigr), \]
where all the cohomology groups are with rational coefficients. Similarly, we get for $q=d-p$,
\[H^{q,q}(\X, \Q) = \ker\Bigl(N\colon H^{q,q}(\X) \to H^{q-1,q+1}(\X)\Bigr) \bigoplus \Im\Bigl(N\colon H^{q+1,q-1}(\X) \to H^{q,q}(\X)\Bigr). \]

Consider as in Section \ref{sec:steenbrink} the polarization $\psi$ on $\ST_1^{\bul, \bul}$ which induces a polarization on the cohomology. The bilinear pairing $\psi$ restricted to $H^{p,p}(\X) \times H^{q,q}(\X)$ coincides moreover with the Poincar\'e duality pairing. The operator $N$ verifies $N^\dual = -N$, for $N^\dual$ the adjoint of $N$ with respect to $\psi$. It follows that the above decompositions of $H^{p,p}(\X)$ and $H^{q,q}(\X, \Q)$ are orthogonal to each other, namely that
\[\ker\Bigl(N\colon H^{p,p}(\X) \to H^{p-1,p+1}(\X)\Bigr) \, \perp \, \Im\Bigl(N\colon H^{q+1,q-1}(\X) \to H^{q,q}(\X)\Bigr), \]
and similarly,
\[\ker\Bigl(N\colon H^{q,q}(\X) \to H^{q-1,q+1}(\X)\Bigr) \, \perp \, \Im\Bigl(N\colon H^{p+1,p-1}(\X) \to H^{p,p}(\X)\Bigr). \]
From the non-degeneracy of $\psi$ on tropical cohomology, we infer that the induced pairing by the polarization between $\ker\Bigl(N\colon H^{p,p}(\X) \to H^{p-1,p+1}(\X)\Bigr)$ and $\ker\Bigl(N\colon H^{q,q}(\X) \to H^{q-1,q+1}(\X)\Bigr)$ is a duality pairing.

By Theorem~\ref{thm:HC} the two kernels above are generated by the classes of tropical cycles of codimension $p$ and $q$, respectively. Moreover, the intersection between cycles is compatible with the polarization evaluated at the pair consisting of the cohomological classes associated to the two cycles. With this preparation, we can now finish the proof.

\begin{proof}[Proof of Theorem~\ref{thm:standard}] Let $p$ be a fixed integer between $0$ and $d$. A tropical cycle of codimension $p$ which is homologically trivial has a trivial intersection with any tropical cycle of codimension $q=d-p$. This proves the implication (Homological equivalence) $\Rightarrow$ (Numerical equivalence).

We now prove the other direction. Consider a tropical cycle $C$ of codimension $p$ in $\X$ and suppose that $C$ has trivial intersection with any tropical cycle of codimension $q=d-p$. It follows from Theorem~\ref{thm:HC} that $\class(C)$ has a trivial pairing with any element of $\ker\bigl(N\colon H^{q,q}(\X) \to H^{q-1,q+1}(\X)\bigr)$. Moreover, it belongs to $\ker\bigl(N\colon H^{p,p}(\X) \to H^{p-1,p+1}(\X)\bigr)$. Since the polarization induces a non-degenerate pairing between the two kernels, we finally infer that $\class(C)$ vanishes. This proves the implication (Numerical equivalence) $\Rightarrow$ (Homological equivalence), and the theorem follows.
\end{proof}

\bibliographystyle{alpha}
\bibliography{$HOME/bibliography/bibliography}

\end{document}